\documentclass[11pt]{article}

\usepackage[margin=1.15in]{geometry}
\usepackage{amsmath,amssymb,amsthm,mathtools}
\usepackage{microtype}
\usepackage{aliascnt}
\usepackage[hidelinks]{hyperref}
\usepackage[nameinlink,capitalize,noabbrev]{cleveref}
\hypersetup{pdftitle={Properties of the Sawin-Wood Model for Global Fields}}

\numberwithin{equation}{section}
\newtheorem{theorem}{Theorem}[section]
\newaliascnt{proposition}{theorem}
\newtheorem{proposition}[proposition]{Proposition}
\aliascntresetthe{proposition}
\newaliascnt{lemma}{theorem}
\newtheorem{lemma}[lemma]{Lemma}
\aliascntresetthe{lemma}
\newaliascnt{corollary}{theorem}

\aliascntresetthe{corollary}
\newaliascnt{question}{theorem}
\newtheorem{question}[question]{Question}
\aliascntresetthe{question}
\theoremstyle{definition}
\newaliascnt{definition}{theorem}
\newtheorem{definition}[definition]{Definition}
\aliascntresetthe{definition}
\theoremstyle{remark}
\newaliascnt{remark}{theorem}

\aliascntresetthe{remark}
\crefname{proposition}{Proposition}{Propositions}
\crefname{lemma}{Lemma}{Lemmas}
\crefname{corollary}{Corollary}{Corollaries}
\crefname{question}{Question}{Questions}
\crefname{definition}{Definition}{Definitions}
\crefname{remark}{Remark}{Remarks}

\DeclareMathOperator{\Aut}{Aut}
\DeclareMathOperator{\End}{End}
\DeclareMathOperator{\Ext}{Ext}
\DeclareMathOperator{\Hom}{Hom}

\DeclareMathOperator{\Sur}{Sur}

\DeclareMathOperator{\GL}{GL}
\newcommand{\F}{\mathbb F}
\newcommand{\Z}{\mathbb Z}
\newcommand{\Prob}{\mathbf P}
\newcommand{\Exp}{\mathbf E}
\newcommand{\GammaAut}{\Aut_{\Gamma}}
\newcommand{\GammaSur}{\Sur_{\Gamma}}
\newcommand{\orSur}{\Sur_{\mathrm{or}}}
\newcommand{\rt}{\rtimes}

\title{Properties of the Sawin--Wood Model for Unramified Galois Groups of Global Fields}
\author{Mark Shusterman \and Gonen Zimmerman}
\date{}

\DeclareMathOperator{\Sym}{Sym}
\newcommand{\ab}{\mathrm{ab}}
\newcommand{\cM}{\mathcal M}
\newcommand{\cI}{\mathcal I}
\newcommand{\eps}{\varepsilon}
\newcommand{\normal}{\trianglelefteq}
\newcommand{\onto}{\twoheadrightarrow}
\newcommand{\qbinom}[2]{\genfrac{[}{]}{0pt}{}{#1}{#2}}

\begin{document}
\maketitle

\begin{abstract}
We prove almost-sure finite generation and finite presentation
in the pro-$|\Gamma|'$ $\Gamma$-category for the Sawin--Wood measure
$\nu_{\Gamma,n,U}$, where $\Gamma$ is a non-trivial finite group, $n$
is coprime to $|\Gamma|$, and $U$ is a nonempty multiset of conjugacy
classes of $\Gamma$.
The $\Gamma$-generator rank of the random group $X$
has the sharp tail
\[
 \Prob(d_\Gamma(X)\ge r)
   =C e^{-\kappa r^2-\zeta r}\bigl(1+O(e^{-\eps r})\bigr),
\]
where $C>0$, $\kappa>0$, and $\zeta\in\mathbb R$ are explicit
constants depending on $n,\Gamma$ and $U$, and $\eps>0$.
Almost surely, $X$ admits a presentation with the minimum number of
$\Gamma$-generators, using at most $|U|$ more defining
$\Gamma$-relations than generators.
We also prove that every open subgroup of $X$ has finite abelianization
almost surely under the assumption that for every complex representation of $\Gamma$ there is an element in $U$ that fixes a nonzero vector. Finally, we prove that
$H_2(X,\widehat{\mathbb Z})=0$ almost surely when
$n=1$, $\Gamma$ is cyclic, and $U$ contains only the conjugacy class of a generator of $\Gamma$ (with multiplicity one).
\end{abstract}

\setcounter{tocdepth}{2}
\tableofcontents

\section{Introduction and basic notation}

The Sawin–Wood measure $\nu_{\Gamma,n,U}$ models the Galois group of the maximal unramified pro-$|\Gamma|'$ extension of a Galois extension $K/\mathbb{Q}$ with $\operatorname{Gal}(K/\mathbb{Q})\cong\Gamma$, as $K$ varies. Here $n=\frac{2}{\gcd(|\Gamma|,2)}$, and $U$ consists of the conjugacy class of complex conjugation in $\Gamma$ (the identity class if $K$ is totally real). For more details, see \cite{SawinWood}. We build on the work of Pagano and Shusterman
\cite{PaganoShusterman}, who use crown theory to establish almost-sure
finite generation for Sawin--Wood's model, and other related random profinite-group models.

The motivation for the \emph{finite-generation and finite-presentation}
results comes, in particular, from a longstanding question of
Shafarevich \cite{Shafarevich}. For a number field $K$ and a finite
set $S$ of finite places, let $K_S/K$ be the maximal Galois extension
unramified at every finite place outside $S$, and put
$G_{K,S}=\operatorname{Gal}(K_S/K)$. Is $G_{K,S}$ topologically
finitely generated? The case $S=\varnothing$, in which no finite
place is allowed to ramify, is especially relevant here.

The function-field analogue of Shafarevich's question is known
to have a positive answer in the unramified case.
Combining the results of Shusterman
\cite{ShustermanBalanced} and Liu
\cite[Theorem~1.2]{LiuPresentations}, we see that for any
global field $K$,
finite topological generation of $G_{K,\varnothing}$ implies
finite topological presentation.
Thus finite generation and finite presentation are the two 
arithmetic finiteness questions to which our almost-sure results 
provide probabilistic analogues.

Fix a finite group $\Gamma$, a positive integer $n$ with
$(n,|\Gamma|)=1$, and a nonempty finite multiset $U$ of conjugacy
classes of $\Gamma$. If $\Gamma=1$, the measure $\nu_{\Gamma,n,U}$
is concentrated on the trivial group (see
\cref{subsec:measurability-trivial}). In this case, the finiteness
and homology-vanishing conclusions are immediate. We henceforth
assume $\Gamma\ne1$.

In this paper, we impose the existence hypothesis of Sawin and Wood in
\cite[Theorem~1.8]{SawinWood}: for each prime $p\mid n$, and for every
nonzero finite-dimensional $\F_p[\Gamma]$-module $V$ there exist $v\in V -\{0\}$ and $\gamma_E\in E\in U$
such that $v^{\gamma_E}=v$.
This fixed-vector condition is vacuous when $n=1$, but we still
assume $U\ne\varnothing$.

A $\Gamma$-group is a group with an action of $\Gamma$ by automorphisms;
in the profinite case the action is continuous. All finite
$\Gamma$-groups appearing in the Sawin--Wood model have
order coprime to $|\Gamma|$, and all profinite $\Gamma$-groups under
consideration are inverse limits of such finite $\Gamma$-groups.  A profinite group is
\emph{small} if it has only
finitely many open subgroups of each index. If $\gamma\in
\Gamma$, write $G^\gamma$ for the fixed-point subgroup.
Choose a representative $\gamma_E$ of each conjugacy class $E$ occurring in $U$,
retaining multiplicities, and define the fixed-point weight
\[
 w_U(G)
   :=\frac{\prod_{E\in U}|G^{\gamma_E}|}{|G^{\Gamma}|}.
\]
This is a positive rational number, independent of the chosen
representatives; the product counts multiplicities in $U$. It is the
group factor in \cite[Theorem~1.8]{SawinWood}.
An $n$-orientation on $G$ is an element of
$H_3(G,\Z/n)^\Gamma$. Although the results stated below concern the
underlying group, this orientation enters the defining moments and
the proof of the second-homology result.

For an oriented profinite $\Gamma$-group $\boldsymbol X=(X,s_X)$
and a finite oriented $\Gamma$-group $\boldsymbol H=(H,s_H)$, write
$\orSur(\boldsymbol X,\boldsymbol H)$ for the set of continuous
$\Gamma$-equivariant epimorphisms $f:X\twoheadrightarrow H$ satisfying
$f_*s_X=s_H$. Unless otherwise stated, $\Hom$ and $\Sur$ denote the
sets of continuous homomorphisms and continuous surjective
homomorphisms, respectively, in the profinite category, including
their equivariant versions. The Sawin--Wood measure $\nu_{\Gamma,n,U}$ is the unique
probability measure on isomorphism classes of small profinite
$n$-oriented pro-$|\Gamma|'$ $\Gamma$-groups characterized by the
following moments: for every finite $n$-oriented $\Gamma$-group
$\boldsymbol H=(H,s_H)$ with $(|H|,|\Gamma|)=1$,
\[
 \Exp_{\boldsymbol X\sim\nu_{\Gamma,n,U}}
   \bigl|\orSur(\boldsymbol X,\boldsymbol H)\bigr|
 =
 \begin{cases}
 \displaystyle
 \frac{|H^2(H\rt\Gamma,\Z/n)|}
      {|H^3(H\rt\Gamma,\Z/n)|\,w_U(H)},
       &H_\Gamma=1,\\[7pt]
 0,    &H_\Gamma\ne1.
 \end{cases}
\]
Here $\Z/n$ has trivial $H\rt\Gamma$-action. Existence and uniqueness
under the above hypotheses follow from
\cite[Theorem~1.8]{SawinWood}. Throughout, $\Exp$ and $\Prob$ denote
expectation and probability with respect to this measure.

For a finite or profinite $\Gamma$-group $G$, let $d_\Gamma(G)$ be the
least cardinality of a set whose $\Gamma$-orbits generate $G$ (topologically
in the profinite case), and set $d_\Gamma(G)=\infty$ if no finite such set
exists. Unless the orientation is explicitly mentioned, the random
variable in the statements below is the underlying $\Gamma$-group
of a sample $\boldsymbol X=(X,s_X)$.

All presentation ranks below are taken in the category of
pro-$|\Gamma|'$ $\Gamma$-groups.

Let $F_m'(\Gamma)$ be the free pro-$|\Gamma|'$ $\Gamma$-group on $m$
$\Gamma$-generators. Thus it is the maximal pro-$|\Gamma|'$ quotient
of the free profinite group on symbols $x_{i,\gamma}$, where
$1\le i\le m$ and $\gamma\in\Gamma$, with
$\sigma(x_{i,\gamma})=x_{i,\sigma\gamma}$. For a finitely
$\Gamma$-generated pro-$|\Gamma|'$ group $H$ and an integer
$m\ge d_\Gamma(H)$, choose a $\Gamma$-equivariant epimorphism
\[
 \pi:F_m'(\Gamma)\twoheadrightarrow H,
\]
and consider the least number of elements generating $\ker\pi$ as a closed
normal $\Gamma$-subgroup of $F_m'(\Gamma)$. For fixed $H$ and $m$, this
number is independent of the choice of $\pi$ by the relation-module
calculation in \cite[Proposition~3.4]{LiuPresentations} and the
normal-generation argument of \cite[Theorem~3.1]{LubotzkyPresentations},
applied with ambient action $F_m'(\Gamma)\rt\Gamma$
(the trivial-kernel case follows because every continuous surjective
endomorphism of a finitely generated profinite group is an automorphism).
Denote it by $r_{\Gamma,m}(H)$.
If $H_\Gamma=1$, then $r_{\Gamma,m}(H)-m$ is independent
of $m\ge d_\Gamma(H)$. Write
\[
 r_\Gamma(H):=r_{\Gamma,d_\Gamma(H)}(H).
\]
With the generators-minus-relations convention, define the relative
deficiency by
\[
 \operatorname{def}_\Gamma(H):=d_\Gamma(H)-r_\Gamma(H).
\]
For $H_\Gamma=1$, this also equals
$\sup_{m\ge d_\Gamma(H)}(m-r_{\Gamma,m}(H))$.

\begin{theorem}\label{thm:sharp-tail}
Let $\cI$ be the set of isomorphism classes of
nontrivial irreducible $\F_p[\Gamma]$-modules, as $p$ ranges over the
primes not dividing $|\Gamma|$. For $V\in\cI$, set
\[
 q_V=|\End_\Gamma(V)|,\qquad
 \tau_V=\dim_{\End_\Gamma(V)}V,\qquad
 u_V=\sum_{E\in U}\dim_{\End_\Gamma(V)}V^{\gamma_E}.
\]
Write $V^\vee=\Hom_{\F_p}(V,\F_p)$ for the dual $\Gamma$-module, and put
\[
 \eta(q)=\prod_{j=1}^\infty(1-q^{-j}).
\]
For $p\nmid n$ or $V\not\simeq V^\vee$, set
\[
 a_V=1,\qquad b_V=u_V,\qquad c_V=\eta(q_V)^{-1}.
\]
For $p\mid n$ and $V\simeq V^\vee$, instead set
\[
 \begin{gathered}
 \lambda_V=\log_{q_V}\left|\Hom_\Gamma\!\left(\bigwedge\nolimits^2_{\F_p}V,\F_p\right)\right|,\\
 a_V=\tfrac12,\qquad b_V=u_V+\tfrac12-\lambda_V,\qquad
 c_V=\left[\eta(q_V)
       \prod_{j=0}^\infty(1+q_V^{-j-u_V-1+\lambda_V})\right]^{-1}.
 \end{gathered}
\]
Now, in both cases, define
\[
 \kappa_V=a_V\tau_V^2\log q_V,\qquad
 \zeta_V=\tau_V\bigl(b_V-2a_V(\tau_V-1)\bigr)\log q_V,
\]
and let
\[
 \kappa=\min_{V\in\cI}\kappa_V,\qquad
 \cM=\{V\in\cI:\kappa_V=\kappa\},\qquad
 \zeta=\min_{V\in\cM}\zeta_V,
\]
\[
 C=\sum_{\substack{V\in\cM\\\zeta_V=\zeta}}
 c_Vq_V^{\,b_V(\tau_V-1)-a_V(\tau_V-1)^2}.
\]

Let $X$ have distribution $\nu_{\Gamma,n,U}$. Then $\kappa>0$,
$\zeta\in\mathbb R$, and $C>0$; moreover, there exists
$\eps>0$, depending only on $\Gamma,n,U$, such that, as the
integer $r$ tends to infinity,
\begin{equation}\label{eq:main-answer}
 \Prob\bigl(d_\Gamma(X)\ge r\bigr)
 =C e^{-\kappa r^2-\zeta r}
   \bigl(1+O_{\Gamma,n,U}(e^{-\eps r})\bigr).
\end{equation}
The same asymptotic holds for $\Prob(d_\Gamma(X)=r)$.
In particular, $X$ is topologically finitely generated as a
$\Gamma$-group with probability one.
\end{theorem}

\begin{theorem}\label{thm:bounded-deficiency}
Let $X$ have distribution $\nu_{\Gamma,n,U}$.
With probability one, $d_\Gamma(X)$ is finite and
\begin{equation}\label{eq:def-bound}
 r_\Gamma(X)\le d_\Gamma(X)+|U|.
\end{equation}
In particular, $X$ is finitely presented as a pro-$|\Gamma|'$
$\Gamma$-group, and
\[
 \operatorname{def}_\Gamma(X)\ge-|U|.
\]
\end{theorem}

The proof of \cref{thm:bounded-deficiency} is given in
\cref{sec:def-presentation}. The finite presentation result uses
the bounds on the deficiency in terms of the second cohomology
proved by Liu and Lubotzky in
\cite{LiuPresentations,LubotzkyPresentations}.

The motivation for the finite abelianizations of open subgroups
comes from the fact that, for any number field $K$ and any finite extension $L/K$, such that $L\leq K^{ur}$, the abelianization of $\operatorname{Gal}(K^{ur}/L)$
is finite by class field theory and Minkowski's theorem.

\begin{theorem}\label{thm:fab}
Let $X$ have distribution $\nu_{\Gamma,n,U}$. Suppose that every
nontrivial irreducible complex $\Gamma$-module $W$ satisfies
\[
 \sum_{E\in U}\dim_{\mathbb C}W^{\gamma_E}>0.
\]
Then, with probability one, every open subgroup $H\le X$ has
finite abelianization.
\end{theorem}

For the next result, $H_2(X,\widehat{\mathbb Z})$ denotes continuous
homology with compact coefficients and trivial $X$-action. 

\begin{theorem}\label{thm:homology-vanishing}
Let $X$ have distribution $\nu_{\Gamma,n,U}$. Suppose that
$n=1$, $\Gamma$ is cyclic, and $U$ contains only the conjugacy class of a generator of $\Gamma$ (with multiplicity one).
Then
\begin{equation}\label{eq:homology-main-vanishing}
 \Prob\bigl(H_2(X,\widehat{\mathbb Z})=0\bigr)=1.
\end{equation}
\end{theorem}

The proofs of \cref{thm:fab,thm:homology-vanishing} are in
\cref{sec:fab,sec:homology}, respectively.
In the number-field case, the hypothesis in
\cref{thm:homology-vanishing} forces
$\Gamma$ to be cyclic of order exactly $2$. Indeed, since the representatives
of $U$ have order at most two, as happens for classes associated
with archimedean places, that hypothesis specializes to
$(\Gamma,n,U)=(C_2,1,\{[\sigma]\})$ with $\sigma\ne1$.

\begin{question}\label{conj:virtual-pronilpotence}
Under the hypothesis of \cref{thm:fab}, is it true that
\[
 \Prob_{\nu_{\Gamma,n,U}}\bigl(X\text{ has an open pronilpotent subgroup}\bigr)=1\,?
\]
\end{question}

Motivated by \cref{conj:virtual-pronilpotence}, one can also ask
the analogous question for number fields.

\begin{question}
\label{conj:arithmetic-virtual-pronilpotence}
Let $K/\mathbb Q$ be a real quadratic extension, and let $L/K$ be the
maximal pro-odd Galois extension unramified at every finite place.
Does $\operatorname{Gal}(L/K)$ necessarily contain a pronilpotent open subgroup?
\end{question}

\subsection{Use of ChatGPT}

GPT6 has been used to assist with coming up with the statements,
conjectures, and proofs.

\subsection{Acknowledgments}

Mark Shusterman is the Dr. A. Edward Friedmann Career Development Chair
in Mathematics.

The authors are supported by the European Union (ERC, Function Fields,
101161909).

The work was supported in part by a grant of access to OpenAI models
through the ChatGPT for Academic Researchers Program.

We thank Carlo Pagano for helpful conversations.

\section{Defining the probability measure}\label{sec:moments}

Sawin and Wood define $\nu_{\Gamma,n,U}$ on small profinite
$n$-oriented $\Gamma$-groups $\boldsymbol X=(X,s_X)$. We use both
this oriented probability law and its pushforward to the underlying
$\Gamma$-group $X$. We begin
by obtaining the underlying moment formula from the oriented one.

\begin{proposition}\label{prop:underlying-moment}
For every finite $\Gamma$-group $H$ whose order is coprime to
$|\Gamma|$,
\begin{equation}
 \Exp\bigl|\GammaSur(X,H)\bigr|
 =
 \begin{cases}
 \displaystyle
 \frac{|H^2(H\rt\Gamma,\Z/n)|}{w_U(H)},
       &H_\Gamma=1,\\[7pt]
 0,    &H_\Gamma\ne1.
 \end{cases}                                      \label{eq:underlying-moment}
\end{equation}
Consequently, the expected number of open normal $\Gamma$-stable subgroups
of $X$ with quotient $\Gamma$-isomorphic to $H$ satisfies
\begin{equation}
 \Exp\#\{N\normal X:N\text{ open and }\Gamma\text{-stable},
                    \,X/N\simeq_\Gamma H\}
 =\frac{\Exp|\GammaSur(X,H)|}{|\GammaAut(H)|}.
\label{eq:quotient-bound}
\end{equation}
This expectation also bounds the probability that such a quotient exists.
There is a set of measure one on which every finite $\Gamma$-quotient
$H$ of $X$ satisfies $H_\Gamma=1$. On the same set, $X_\Gamma=1$.
\end{proposition}

\begin{proof}
Let $\boldsymbol X=(X,s_X)$.  For each underlying $\Gamma$-surjection
$f:X\twoheadrightarrow H$, exactly one orientation on $H$, namely
$f_*s_X$, makes $f$ orientation preserving.  Hence
\[
 |\GammaSur(X,H)|
   =\sum_{s\in H_3(H,\Z/n)^\Gamma}
       |\orSur(\boldsymbol X,(H,s))|.
\]
By \cite[Theorem~1.8]{SawinWood},
\[
 \Exp\bigl|\orSur(\boldsymbol X,(H,s))\bigr| = \mathbf 1_{H_\Gamma=1}
 \frac{|H^2(H\rt\Gamma,\Z/n)|}
 {|H^3(H\rt\Gamma,\Z/n)|\,w_U(H)}.
\]
Thus
\[
\Exp\bigl|\GammaSur(X,H)\bigr|
=
\mathbf 1_{H_\Gamma=1}
\frac{|H_3(H,\Z/n)^\Gamma|\,
      |H^2(H\rt\Gamma,\Z/n)|}
     {|H^3(H\rt\Gamma,\Z/n)|\,w_U(H)}.
\]
By \cite[Lemmas~5.4 and~5.5]{SawinWood} (cf. the proof of
\cite[Lemma~4.7]{SawinWood}),
\[
|H_3(H,\Z/n)^\Gamma|
=
|H^3(H\rt\Gamma,\Z/n)|.
\]
These factors cancel, proving \eqref{eq:underlying-moment}.

Each open normal $\Gamma$-stable subgroup $N$ with $X/N\simeq_\Gamma H$
is the kernel of exactly $|\GammaAut(H)|$ $\Gamma$-equivariant
surjections $X\to H$: these correspond to the $\Gamma$-isomorphisms
$X/N\to H$, which differ by postcomposition with a unique
$\Gamma$-automorphism of $H$. Counting surjections by their kernels,
dividing by $|\GammaAut(H)|$, and taking expectations gives
\eqref{eq:quotient-bound}.

The probability that the quotient count is positive is at most its
expectation, by Markov's inequality.

Finally, there are only countably many finite $\Gamma$-groups up to
isomorphism. The union of the zero-moment quotient events with
$H_\Gamma\ne1$ has measure zero. A nontrivial profinite quotient
$X_\Gamma$ would have a nontrivial finite quotient with trivial
$\Gamma$-action, contradicting this conclusion. Hence $X_\Gamma=1$
on the same full-measure set.
\end{proof}

\subsection{Measurability and the trivial case}\label{subsec:measurability-trivial}

Finite-quotient and generator-rank events considered in the paper are
measurable. The former are finite unions over orientations of finite-level
quotient events in the Sawin--Wood space, and the following lemma
shows that the latter are countable unions of the former.

\begin{lemma}\label{lem:compactness}
For every profinite $\Gamma$-group $Y$ and every nonnegative integer $r$,
\[
 d_\Gamma(Y)\le r \quad\Longleftrightarrow\quad
 d_\Gamma(Y/N)\le r
 \quad\text{for every open normal $\Gamma$-subgroup $N$ of $Y$.}
\]
\end{lemma}

\begin{proof}
Only the reverse implication needs proof.  Suppose every finite
$\Gamma$-quotient is generated by $r$ $\Gamma$-orbits.  For each open
normal $\Gamma$-subgroup $N$, let $\mathcal G_N\subseteq Y^r$ be the set
of tuples whose images $\Gamma$-generate $Y/N$.  It is a nonempty clopen
subset.  Any finite intersection of the $\mathcal G_N$ contains a tuple
which generates the quotient by the intersection of the corresponding
$N$'s.  Compactness of $Y^r$ gives a point in
$\bigcap_N\mathcal G_N$, and that tuple topologically
$\Gamma$-generates $Y$.
\end{proof}

In particular, no finite-generation assumption is implicit in these
arguments.

To justify the claim in the introduction, suppose $\Gamma=1$.
Then $H_\Gamma=H$ for every finite quotient $H$ of $X$. By
\cref{prop:underlying-moment}, all such quotients are trivial almost
surely, and hence $X=1$ almost surely.

\section{Proof of the sharp-tail theorem}\label{sec:sharp-tail-proof}

This section proves \cref{thm:sharp-tail}. We begin with the constants
and the elementary abelian multiplicity laws that enter the argument.

\subsection{Constants in the sharp-tail asymptotic}\label{subsec:tail-parameters}

The constants in \cref{thm:sharp-tail} are determined by the nontrivial
irreducible representations of $\Gamma$. We first calculate their
multiplicity distributions, then combine their leading contributions
to obtain the generator-rank tail.

We use the set $\cI$ and the parameters $q_V,\tau_V,u_V$ defined in
\cref{thm:sharp-tail}. As in the definition of $w_U$, the sum defining
$u_V$ retains the multiplicities in $U$ and does not depend on the
chosen representatives $\gamma_E$. The dual $V^\vee$ has the
contragredient action.

\subsubsection{Calculating \texorpdfstring{$\Prob(V\text{ has multiplicity }k\text{ in }X^{\ab}/pX^{\ab})$}{the probability that V has multiplicity k in the maximal elementary abelian p-quotient of X}}\label{subsubsec:multiplicity-distribution}

This subsection explains where the constants defined in
\cref{thm:sharp-tail} come from. The multiplicity formulas are
proved in \cref{sec:abelian}, and the theorem in \cref{sec:conclusion}.

Let $m_V$ be the multiplicity of $V$ in $X^{\ab}/pX^{\ab}$, the maximal
elementary abelian $p$-quotient of $X$. The summation argument in
\cref{lem:individual-tail} gives
\begin{equation}\label{eq:multiplicity-inversion}
 \Prob(m_V=k)=\sum_{t=k}^{\infty}
 (-1)^{t-k}q_V^{\binom{t-k}{2}}\qbinom{t}{k}_{q_V}
 \frac{|H^2(V^t\rt\Gamma,\Z/n)|}
 {|\GammaAut(V^t)|\,w_U(V^t)}.
\end{equation}
Here $\qbinom{t}{k}_q$ is the number of $k$-dimensional subspaces
of $\F_q^t$. By \cref{prop:underlying-moment}, the ratio is the expected
number of open normal $\Gamma$-stable subgroups of $X$ with quotient
$\Gamma$-isomorphic to $V^t$. The alternating sum extracts the
probability that the multiplicity is exactly $k$.

For the cohomology factor, if $p\mid n$ and $V\simeq V^\vee$,
use $\lambda_V$ as defined in \cref{thm:sharp-tail}. Then
\[
 H^2(V^t\rt\Gamma,\Z/n)\simeq
 \Hom_\Gamma\!\left(\bigwedge\nolimits^2_{\F_p}(V^t),\F_p\right),\qquad
 |H^2(V^t\rt\Gamma,\Z/n)|=q_V^{\binom{t}{2}+\lambda_Vt}.
\]
If $p\nmid n$ or $V\not\simeq V^\vee$, then
\[
 H^2(V^t\rt\Gamma,\Z/n)=0,\qquad
 |H^2(V^t\rt\Gamma,\Z/n)|=1.
\]
The calculation is proved in \cref{lem:multiplicity-moments}.

The weight and automorphism factors are
\[
 w_U(V^t)=q_V^{u_Vt},\qquad
 |\GammaAut(V^t)|=|\GL_t(q_V)|=q_V^{t^2}\eta_{q_V}(t),
\]
and, for a prime power $q$, we set
\begin{equation}\label{eq:eta}
 \eta_q(k)=\prod_{j=1}^{k}(1-q^{-j}),\qquad \eta_q(0)=1.
\end{equation}
Evaluating \eqref{eq:multiplicity-inversion} gives the following formulas
for every $k\ge0$, proved in \cref{lem:individual-tail}.
If $p\nmid n$ or $V\not\simeq V^\vee$, then
\begin{equation}\label{eq:ordinary-mass}
 \Prob(m_V=k)
 =\frac{q_V^{-k^2-u_Vk}}{\eta_{q_V}(k)}
   \prod_{j=1}^{\infty}(1-q_V^{-k-u_V-j}).
\end{equation}
If $p\mid n$ and $V\simeq V^\vee$, then
\begin{equation}\label{eq:exceptional-mass}
 \Prob(m_V=k)
 =\frac{q_V^{-k^2/2-(u_V+1/2-\lambda_V)k}}{\eta_{q_V}(k)}
   \prod_{j=0}^{\infty}(1+q_V^{-j-u_V-1+\lambda_V})^{-1}.
\end{equation}

The constants $a_V,b_V,c_V$ defined in \cref{thm:sharp-tail}
give the quadratic exponent, linear exponent, and leading factor in the
multiplicity tail.
Schur's lemma gives $0\le\lambda_V\le1$.
The Sawin--Wood fixed-vector hypothesis gives $u_V\ge1$ when $p\mid n$, so all the parameters
are well-defined, the products are positive, and $b_V\ge0$.
The tail has the same leading term as its first mass: by
\cref{lem:individual-tail}, as $k\to\infty$,
\[
 \Prob(m_V\ge k)=c_Vq_V^{-a_Vk^2-b_Vk}
                    \bigl(1+O_V(q_V^{-k})\bigr).
\]

\subsubsection{The main terms in the generator rank tail}\label{subsubsec:generator-tail-constants}

For $V\in\cI$ over $\F_p$, the group
$V^{m_V}$ is the $V$-isotypic quotient of $X^{\ab}/pX^{\ab}$.
Since $d_\Gamma(V^{m_V})=\lceil m_V/\tau_V\rceil$, the preceding multiplicity
formula gives, as the integer $r\to\infty$,
\[
 \begin{aligned}
 \Prob\bigl(d_\Gamma(V^{m_V})\ge r\bigr)
   &=\Prob\bigl(m_V\ge \tau_V(r-1)+1\bigr)\\
 &=c_Vq_V^{b_V(\tau_V-1)-a_V(\tau_V-1)^2}
   e^{-\left(a_V\tau_V^2 r^2
      +\tau_V(b_V-2a_V(\tau_V-1))r\right)\log q_V}
   \bigl(1+O_V(q_V^{-r})\bigr).
 \end{aligned}
\]
The quadratic and linear coefficients in the decay exponent are
$\kappa_V$ and $\zeta_V$ from \cref{thm:sharp-tail}.
Here and throughout, $\log$ is the natural logarithm.
The smallest quadratic coefficient determines the leading decay,
as reflected in the definitions of $\kappa$ and $\cM$.
The set $\cI$ is nonempty.  The minimum is attained, $\kappa>0$, and
$\cM$ is finite: since $a_V\ge\tfrac12$, an upper bound on $\kappa_V$
bounds $q_V$ and $\tau_V$, and there are
finitely many irreducibles in any given characteristic.
Among these modules, those with the smallest linear coefficient
contribute to the leading constant $C$ defined in \cref{thm:sharp-tail}.
Write
\begin{equation}\label{eq:threshold}
 k_V(r)=\tau_V(r-1)+1.
\end{equation}
The proof gives the more precise expansion
\begin{equation}\label{eq:finite-sum-answer}
 \Prob\bigl(d_\Gamma(X)\ge r\bigr)
 =\left(\sum_{V\in\cM}
 c_Vq_V^{-a_Vk_V(r)^2-b_Vk_V(r)}\right)
 \bigl(1+O_{\Gamma,n,U}(e^{-\eps r})\bigr).
\end{equation}
\subsection{Calculating the asymptotic for \texorpdfstring{$\Prob(d_\Gamma(X^{\ab})\ge r)$}{the probability that the abelian generator rank is at least r}}\label{sec:abelian}

\Cref{sec:relative,sec:rigidity,sec:estimates} will show that the difference
between $\Prob(d_\Gamma(X)\ge r)$ and $\Prob(d_\Gamma(X^{\ab})\ge r)$
is negligible relative to the abelian tail calculated here.

\begin{lemma}\label{lem:abelian-rank}
For each $V\in\cI$ over $\F_p$, let $m_V$ be the multiplicity of $V$
in the $\F_p[\Gamma]$-module $X^{\ab}/pX^{\ab}$, the maximal
elementary abelian $p$-quotient of $X$. Almost surely,
\begin{equation}\label{eq:abelian-rank}
 d_\Gamma(X^{\ab})=\sup_{V\in\cI}
                  \left\lceil\frac{m_V}{\tau_V}\right\rceil.
\end{equation}
In particular, $d_\Gamma(X^{\ab})\ge r$ if and only if
$m_V\ge k_V(r)$ for at least one $V$.
\end{lemma}
\begin{proof}
For $p\nmid|\Gamma|$, the maximal elementary abelian $p$-quotient of a
small profinite group is finite. Indeed, every nonzero continuous map
to $C_p$ has an open kernel of index $p$, and each such kernel supports
only $p-1$ such maps. Thus $\Hom(X,C_p)$ is a finite-dimensional
$\F_p$-space. On the full-measure set of \cref{prop:underlying-moment},
there are no trivial constituents.

By the Wedderburn decomposition of $\F_p[\Gamma]$, its regular module
contains $\tau_V$ copies of $V$.  A semisimple module with multiplicities
$m_V$ is generated by $k$ elements if and only if $m_V\le k\tau_V$ for all
$V$.  For a finite abelian $p$-group with a $\Gamma$-action, a tuple
generates if and only if it generates modulo $p$: if the subgroup it
$\Gamma$-generates is $N$ and $G=N+pG$, iteration gives $G=N$.
Primary components can be generated simultaneously.  Apply these facts
to finite abelian quotients and use \cref{lem:compactness}.
\end{proof}

For integers $m,t\ge0$, we have
\[
 |\Sur_{\F_q}(\F_q^m,\F_q^t)|
         =\prod_{j=0}^{t-1}(q^m-q^j).
\]
The count is $1$ when $t=0$ and $0$ when $t>m$.

\begin{lemma}\label{lem:multiplicity-moments}
Fix $V\in\cI$. Then
\begin{equation}\label{eq:multiplicity-moments}
 \Exp |\Sur_{\F_{q_V}}(\F_{q_V}^{m_V},\F_{q_V}^{t})|=
 \begin{cases}
  q_V^{-u_Vt},& p\nmid n\text{ or }V\not\simeq V^\vee,\\
  q_V^{\binom t2+(\lambda_V-u_V)t},&p\mid n,
                                      \ V\simeq V^\vee.
 \end{cases}
\end{equation}
\end{lemma}
\begin{proof}
Every continuous $\Gamma$-equivariant map from $X$ to $V^t$ factors
through the $V$-isotypic quotient $V^{m_V}$, so
$|\GammaSur(X,V^t)|=|\Sur_{\F_{q_V}}(\F_{q_V}^{m_V},\F_{q_V}^{t})|$.
Apply \eqref{eq:underlying-moment} to $G=V^t$. Its weight is $q_V^{u_Vt}$
and its coinvariants vanish.  If $p\nmid n$, its positive cohomology
with coefficients $\Z/n$ vanishes.  Otherwise only the $p$-primary part
of the coefficients contributes.  For $p^e\Vert n$, the universal
coefficient sequence is
\[
 0\longrightarrow \Ext^1_\Z(V^t,\Z/p^e)
 \longrightarrow H^2(V^t,\Z/p^e)
 \longrightarrow \Hom_\Z(\bigwedge^2 V^t,\Z/p^e)
 \longrightarrow0.
\]
The left term is a direct sum of copies of $V^\vee$ and has no
$\Gamma$-invariants.  Taking invariants is exact because
$p\nmid|\Gamma|$.  The decomposition
\[
 \bigwedge\nolimits^2 V^t\simeq
 \bigoplus_{i=1}^t\bigwedge\nolimits^2 V
 \ \oplus\!\bigoplus_{i<j}(V\otimes_{\F_p}V)
\]
is valid also for $p=2$. Each cross summand has $q_V$ invariant forms
when $V$ is self-dual and only the zero form otherwise.  The diagonal
summands contribute $q_V^{\lambda_Vt}$ in the self-dual case, and only
the zero invariant form otherwise. Averaging over $\Gamma$
identifies this calculation with the cohomology of $V^t\rt\Gamma$.
\end{proof}

\subsubsection{The multiplicity tails of irreducible \texorpdfstring{$\Gamma$}{Gamma}-modules}

\begin{lemma}\label{lem:individual-tail}
For each fixed $V\in\cI$,
\begin{equation}\label{eq:individual-tail}
 \Prob(m_V\ge k)=c_Vq_V^{-a_Vk^2-b_Vk}
                         \bigl(1+O_V(q_V^{-k})\bigr).
\end{equation}
The implicit constant may also depend on $n,U$.  Uniformly in $V$ and $k\ge0$,
\begin{equation}\label{eq:uniform-multiplicity-bound}
 \Prob(m_V\ge k)\le \rho^{-1}q_V^{-a_Vk^2-b_Vk},
 \qquad \rho=\eta(2)>0.
\end{equation}
\end{lemma}
\begin{proof}
First, $|\Sur_{\F_{q_V}}(\F_{q_V}^{m_V},\F_{q_V}^{k})|\ge|\GL_k(q_V)|$ whenever $m_V\ge k$.
Markov's inequality and \cref{lem:multiplicity-moments} give
\[
 \Prob(m_V\ge k)\le
 \frac{\Exp |\Sur_{\F_{q_V}}(\F_{q_V}^{m_V},\F_{q_V}^{k})|}{|\GL_k(q_V)|}
 =\frac{q_V^{-a_Vk^2-b_Vk}}{\eta_{q_V}(k)}
 \le\rho^{-1}q_V^{-a_Vk^2-b_Vk},
\]
since $\eta_{q_V}(k)\ge\eta(q_V)\ge\rho$. This proves the uniform bound.

To recover the individual masses, use
\[
 \qbinom{m}{t}_{q_V}
 =\frac{|\Sur_{\F_{q_V}}(\F_{q_V}^{m},\F_{q_V}^{t})|}{|\GL_t(q_V)|},
\]
with value zero for $t>m$. The finite $q$-binomial theorem gives
\[
 \mathbf1_{\{m=k\}}=
 \sum_{e\ge0}(-1)^e q_V^{\binom e2}
     \qbinom{m}{k+e}_{q_V}\qbinom{k+e}{k}_{q_V}.
\]
Indeed, for $m\ge k$ the sum equals
$\qbinom{m}{k}_{q_V}\prod_{j=0}^{m-k-1}(1-q_V^j)$;
it is $1$ when $m=k$ and vanishes otherwise. For $m<k$ all terms
are zero. Thus the alternating sum isolates multiplicity exactly $k$.

Taking expectations and applying \eqref{eq:underlying-moment} gives
\eqref{eq:multiplicity-inversion}, which we repeat here:
\[
 \Prob(m_V=k)=\sum_{t=k}^{\infty}
 (-1)^{t-k}q_V^{\binom{t-k}{2}}\qbinom{t}{k}_{q_V}
 \frac{|H^2(V^t\rt\Gamma,\Z/n)|}
 {|\GammaAut(V^t)|\,w_U(V^t)}.
\]
Writing $t=k+e$ and using the formula for the Gaussian binomial
coefficient gives
\begin{equation}\label{eq:moment-inversion}
 \Prob(m_V=k)=\frac1{\eta_{q_V}(k)}
 \sum_{e\ge0}
 \frac{(-1)^e q_V^{-k^2-ke-e(e+1)/2}}{\eta_{q_V}(e)}
       \Exp |\Sur_{\F_{q_V}}(\F_{q_V}^{m_V},\F_{q_V}^{k+e})|.
\end{equation}
Substituting the moments from \cref{lem:multiplicity-moments} and
separating the terms involving $k$ from those involving $e$ gives
\[
 \Prob(m_V=k)=
 \begin{cases}
 \displaystyle
 \frac{q_V^{-k^2-u_Vk}}{\eta_{q_V}(k)}
 \sum_{e\ge0}\frac{(-1)^e q_V^{-(k+u_V)e-e(e+1)/2}}
                         {\eta_{q_V}(e)},
       &p\nmid n\text{ or }V\not\simeq V^\vee,\\[1.2em]
 \displaystyle
 \frac{q_V^{-k^2/2-(u_V+1/2-\lambda_V)k}}{\eta_{q_V}(k)}
 \sum_{e\ge0}\frac{(-1)^e q_V^{-(u_V+1-\lambda_V)e}}
                         {\eta_{q_V}(e)},
       &p\mid n,\ V\simeq V^\vee.
 \end{cases}
\]
These expressions also justify the
interchange of sum and expectation: the first series converges
absolutely by its negative quadratic exponent in $e$, and the
absolute terms of the second are bounded by a geometric series,
since $u_V+1-\lambda_V\ge1$ and $\eta_{q_V}(e)\ge\rho$.

Euler's product identities \cite[Corollary~2.2, p.~19]{AndrewsPartitions} are
\[
 \sum_{e\ge0}\frac{(-1)^e q_V^{-e(e+1)/2}z^e}{\eta_{q_V}(e)}
 =\prod_{j\ge1}(1-zq_V^{-j}),\quad z\in\mathbb C,\qquad
 \sum_{e\ge0}\frac{z^e}{\eta_{q_V}(e)}
 =\prod_{j\ge0}(1-zq_V^{-j})^{-1},\quad |z|<1.
\]
In the first case, substitute $z=q_V^{-k-u_V}$ into the first identity;
in the second, substitute $z=-q_V^{-u_V-1+\lambda_V}$ into the second.
Equation \eqref{eq:moment-inversion} then gives exactly
\eqref{eq:ordinary-mass} and \eqref{eq:exceptional-mass}.

Since $\eta_{q_V}(k)=\eta(q_V)(1+O(q_V^{-k}))$ and
$\prod_{j\ge1}(1-q_V^{-k-u_V-j})=1+O(q_V^{-k})$,
these mass formulas give
\[
 \Prob(m_V=k)=c_Vq_V^{-a_Vk^2-b_Vk}
                  \bigl(1+O_V(q_V^{-k})\bigr).
\]
Finally, the uniform bound at $k+1$ gives
\[
 0\le\Prob(m_V\ge k)-\Prob(m_V=k)
 =\Prob(m_V\ge k+1)
 \le\rho^{-1}q_V^{-a_V(k+1)^2-b_V(k+1)}.
\]
Since $2a_V\ge1$ and $b_V\ge0$, this is a relative $O_V(q_V^{-k})$
error. This proves \eqref{eq:individual-tail}.
\end{proof}

\subsubsection{The asymptotic for \texorpdfstring{$\Prob(d_\Gamma(X^{\ab})\ge r)$}{the probability that the abelian generator rank is at least r}}

\begin{proposition}\label{prop:abelian-tail}
The right side of \eqref{eq:finite-sum-answer} is an asymptotic for
$\Prob(d_\Gamma(X^{\ab})\ge r)$, with the same type of relative error.
\end{proposition}
\begin{proof}
By \cref{lem:abelian-rank},
\[
 \Prob\bigl(d_\Gamma(X^{\ab})\ge r\bigr)
 =\Prob\left(\bigcup_{V\in\cI}\{m_V\ge k_V(r)\}\right).
\]
The union bound and the first two terms of inclusion--exclusion give
\begin{equation}\label{eq:abelian-union-error}
 \begin{aligned}
 &\left|\Prob\bigl(d_\Gamma(X^{\ab})\ge r\bigr)
      -\sum_{V\in\cM}\Prob\bigl(m_V\ge k_V(r)\bigr)\right|\le\\
 &
   \sum_{\substack{\{V,W\}\subseteq\cM\\V\ne W}}
     \Prob\bigl(m_V\ge k_V(r),\ m_W\ge k_W(r)\bigr)+
   \sum_{V\in\cI\setminus\cM}\Prob\bigl(m_V\ge k_V(r)\bigr).
 \end{aligned}
\end{equation}
The first sum on the right bounds the overcounting when we add the
individual probabilities for $V\in\cM$, while the second bounds the
contribution from modules outside $\cM$. All pair sums below are over
unordered pairs of distinct modules. We split the first sum as
\[
 \begin{aligned}
 &\sum_{\substack{\{V,W\}\subseteq\cM\\V\ne W}}
   \Prob\bigl(m_V\ge k_V(r),\ m_W\ge k_W(r)\bigr)=\\
 &
   \sum_{\substack{\{V,W\}\subseteq\cM,\ V\ne W\\
          W\not\simeq V^\vee\ \text{or}\ \gcd(|V|,n)=1}}
   \Prob\bigl(m_V\ge k_V(r),\ m_W\ge k_W(r)\bigr)+
   \sum_{\substack{\{V,W\}\subseteq\cM,\ V\ne W\\
          W\simeq V^\vee,\ \gcd(|V|,n)>1}}
   \Prob\bigl(m_V\ge k_V(r),\ m_W\ge k_W(r)\bigr).
 \end{aligned}
\]

By \cref{lem:abelian-rank}, the threshold is
$k_V(r)=\tau_V(r-1)+1$.
Since $b_V\ge0$ and $\kappa_V=a_V\tau_V^2\log q_V$, we will use
\[
 q_V^{-a_Vk_V(r)^2-b_Vk_V(r)}
 \le q_V^{-a_V\tau_V^2(r-1)^2}
 =e^{-\kappa_V(r-1)^2}.
\]
For comparison, expanding $k_V(r)=\tau_Vr-(\tau_V-1)$ gives the main term
\[
 \begin{aligned}
 \sum_{V\in\cM}c_Vq_V^{-a_Vk_V(r)^2-b_Vk_V(r)} &= e^{-\kappa r^2-\zeta r}
   \sum_{V\in\cM}
     c_Vq_V^{b_V(\tau_V-1)-a_V(\tau_V-1)^2}e^{-(\zeta_V-\zeta)r}\\
 &\asymp_{\Gamma,n,U}e^{-\kappa r^2-\zeta r}.
 \end{aligned}
\]
The last comparison holds because $\cM$ is finite, every coefficient
is positive, and $\zeta=\min_{V\in\cM}\zeta_V$.

We first bound the sum over pairs with $W\not\simeq V^\vee$ or
$\gcd(|V|,n)=1$. For these pairs, the moment of the target
$V^k\oplus W^l$ factors into the product of the individual moments.
This follows from the same exterior-power calculation as in
\cref{lem:multiplicity-moments}; for different primes there are no
cross terms. Markov's inequality gives
\begin{equation}\label{eq:independent-intersection}
 \begin{aligned}
 \Prob(m_V\ge k,m_W\ge l)
 &\le\frac{\Exp\bigl(|\Sur_{\F_{q_V}}(\F_{q_V}^{m_V},\F_{q_V}^{k})|\,|\Sur_{\F_{q_W}}(\F_{q_W}^{m_W},\F_{q_W}^{l})|\bigr)}
               {|\GL_k(q_V)|\,|\GL_l(q_W)|}\\
 &=\frac{q_V^{-a_Vk^2-b_Vk}q_W^{-a_Wl^2-b_Wl}}
              {\eta_{q_V}(k)\eta_{q_W}(l)}\le\rho^{-2}q_V^{-a_Vk^2-b_Vk}q_W^{-a_Wl^2-b_Wl}.
 \end{aligned}
\end{equation}
Applying the preceding exponential bound to both factors, and using
$\kappa_V=\kappa_W=\kappa$ for $V,W\in\cM$, gives
\[
 \begin{aligned}
 &\sum_{\substack{\{V,W\}\subseteq\cM,\ V\ne W\\
          W\not\simeq V^\vee\ \text{or}\ \gcd(|V|,n)=1}}
   \Prob\bigl(m_V\ge k_V(r),\ m_W\ge k_W(r)\bigr)\\
 &\quad\le\rho^{-2}
   \sum_{\substack{\{V,W\}\subseteq\cM,\ V\ne W\\
          W\not\simeq V^\vee\ \text{or}\ \gcd(|V|,n)=1}}
       e^{-(\kappa_V+\kappa_W)(r-1)^2}\\
 &\quad\le\rho^{-2}\binom{|\cM|}{2}
                    e^{-2\kappa(r-1)^2}.
 \end{aligned}
\]
There are at most $\binom{|\cM|}{2}$ pairs. Relative to the main term,
this bound has size $O(e^{-\kappa r^2/2})$, since
\[
 \frac{e^{-2\kappa(r-1)^2}}{e^{-\kappa r^2-\zeta r}}
 =e^{-\kappa r^2+(4\kappa+\zeta)r-2\kappa}
 =O\bigl(e^{-\kappa r^2/2}\bigr).
\]
No independence of the random variables is required.

We now bound the sum over pairs with $W\simeq V^\vee$ and
$\gcd(|V|,n)>1$. These are dual pairs in a characteristic $p\mid n$.
Here $q_W=q_V$, $\tau_W=\tau_V$, and $u_W=u_V$, so
$k_W(r)=k_V(r)$. The mixed moments are
\[
 \Exp\bigl(|\Sur_{\F_{q_V}}(\F_{q_V}^{m_V},\F_{q_V}^{k})|\,|\Sur_{\F_{q_V}}(\F_{q_V}^{m_W},\F_{q_V}^{l})|\bigr)
 =q_V^{kl-u_V(k+l)}.
\]
Thus
\begin{equation}\label{eq:dual-intersection}
 \begin{aligned}
 \Prob\bigl(m_V\ge k_V(r),\ m_W\ge k_W(r)\bigr)\le&
 \frac{\Exp\bigl(|\Sur_{\F_{q_V}}(\F_{q_V}^{m_V},\F_{q_V}^{k_V(r)})|
                  \,|\Sur_{\F_{q_V}}(\F_{q_V}^{m_W},\F_{q_V}^{k_V(r)})|\bigr)}
      {|\GL_{k_V(r)}(q_V)|^2}
      \\&=\frac{q_V^{-k_V(r)^2-2u_Vk_V(r)}}
                   {\eta_{q_V}(k_V(r))^2}
 \le\rho^{-2}q_V^{-u_Vk_V(r)}
                      q_V^{-k_V(r)^2-u_Vk_V(r)}.
 \end{aligned}
\end{equation}
Here $u_V\ge1$, $a_V=1$, $b_V=u_V$, and
$c_V=\eta(q_V)^{-1}\ge1$. Moreover,
\[
 k_V(r)-r=(\tau_V-1)(r-1)\ge0,
 \qquad
 q_V^{-u_Vk_V(r)}\le2^{-k_V(r)}\le2^{-r}.
\]
Each dual pair has two members with the same parameters. Replacing
the sum over pairs by the sum over their members therefore gives
\[
 \begin{aligned}
 \sum_{\substack{\{V,W\}\subseteq\cM,\ V\ne W\\
          W\simeq V^\vee,\ \gcd(|V|,n)>1}}
   \Prob\bigl(m_V\ge k_V(r),\ m_W\ge k_W(r)\bigr)
 &\le\rho^{-2}
   \sum_{\substack{V\in\cM,\ V^\vee\in\cM\\
                   V\not\simeq V^\vee,\ \gcd(|V|,n)>1}}
       q_V^{-u_Vk_V(r)}q_V^{-k_V(r)^2-u_Vk_V(r)}\\
 &\le\rho^{-2}2^{-r}
       \sum_{V\in\cM}c_Vq_V^{-a_Vk_V(r)^2-b_Vk_V(r)}.
 \end{aligned}
\]
Thus this error is at most $\rho^{-2}2^{-r}$ times the main term.
In particular, both members of a dual pair are to be included in the
leading sum.

Finally, we bound the sum over $V\notin\cM$ in
\eqref{eq:abelian-union-error}.
The set of quadratic decay coefficients $\{\kappa_V:V\in\cI\}$ is locally finite. Let
$\kappa_2>\kappa$ be its least value greater than $\kappa$.  It exists
because there are nontrivial modules in arbitrarily large
characteristics. By \eqref{eq:uniform-multiplicity-bound} and the
exponential bound above, for $r\ge3$ we have
\[
 \begin{aligned}
 \sum_{V\notin\cM}\Prob\bigl(m_V\ge k_V(r)\bigr)
 &\le\rho^{-1}\sum_{V\notin\cM}e^{-\kappa_V(r-1)^2}\\
 &=\rho^{-1}\sum_{V\notin\cM}
       e^{-\kappa_V((r-1)^2-4)}e^{-4\kappa_V}\\
 &\le\rho^{-1}e^{-\kappa_2((r-1)^2-4)}
       \sum_{V\in\cI}e^{-4\kappa_V}\\
 &\ll_{\Gamma,n}e^{-\kappa_2(r-1)^2}.
 \end{aligned}
\]
In the third line we used $\kappa_V\ge\kappa_2$ for $V\notin\cM$
and $(r-1)^2-4\ge0$. The remaining sum is finite: there are at most
$|\Gamma|$ simple modules in each characteristic, and
$\kappa_V\ge\tfrac12\log p$, so
\[
 \sum_{V\in\cI}e^{-4\kappa_V}
 \le |\Gamma|\sum_{p\nmid|\Gamma|}p^{-2}<\infty.
\]
Finally, $\kappa_2>\kappa$ makes this error smaller than the main term:
\[
 \frac{e^{-\kappa_2(r-1)^2}}{e^{-\kappa r^2-\zeta r}}
 =e^{-(\kappa_2-\kappa)r^2+(2\kappa_2+\zeta)r-\kappa_2}
 =O\bigl(e^{-(\kappa_2-\kappa)r^2/2}\bigr).
\]

Apply \eqref{eq:abelian-union-error} and the
individual asymptotics \eqref{eq:individual-tail}. Since $\cM$ is finite
and all thresholds grow at least linearly with $r$, the resulting
relative error is $O(e^{-\eps r})$ for some $\eps>0$.
\end{proof}

\subsection{Producing finite quotient groups for \texorpdfstring{$d_\Gamma(X)\ge r$}{generator rank at least r}}\label{sec:relative}

For $r\ge4$, the event $d_\Gamma(X)\ge r$ forces a
$\Gamma$-equivariant epimorphism from $X$ onto one of the finite groups
$\pi_j(G)$ constructed in this subsection, on the full-measure set of
\cref{prop:underlying-moment}. By \cref{lem:compactness}, we first choose a
finite $\Gamma$-quotient $G$ with $d_\Gamma(G)\ge r$, then apply crown
reduction and project to a prescribed multiplicity.

The projected groups $\pi_j(G)$ need not satisfy
$d_\Gamma(\pi_j(G))\ge r$. We retain the implication that a group $X$
with $d_\Gamma(X)\ge r$ has one of these groups as a quotient.
We will use the quotient moments from \cref{prop:underlying-moment}
to bound the probability that these quotient groups occur.

The finite-group and cohomological inputs below rely on the
classification of finite simple groups. Consequently, our proofs of
\cref{thm:sharp-tail,thm:bounded-deficiency} also rely on it.

\subsubsection{Crown theory with \texorpdfstring{$\Gamma$}{Gamma}-action}\label{sec:crown-gamma}

We state the crown arguments relative to an arbitrary subgroup.
For a finite group $E$ and a subgroup $H\le E$, put
$d_H(E)=\min\{k:E=\langle H,x_1,\ldots,x_k\rangle\}$.
Such a tuple is called a \emph{relative generating tuple}.
For $E=G\rt\Gamma$ with its specified $H=\Gamma$, this equals
$d_\Gamma(G)$. The results below require no coprimality assumption.

\begin{lemma}\label{lem:relative-lifting}
If $N\normal E$ and $k\ge d_H(E)$, every relative generating $k$-tuple
of $E/N$ lifts to one of $E$.
\end{lemma}
\begin{proof}
For subgroups $M\le K\le E$, let $\mu_E(M,K)$ denote the M\"obius
function of the subgroup lattice of $E$, characterized by
$\mu_E(M,M)=1$ and $\sum_{M\le J\le K}\mu_E(M,J)=0$ whenever $M<K$.
M\"obius inversion on the subgroup lattice counts its generating lifts as
\[
 \sum_{\substack{M\le E,\ H\le M\\MN=E}}
                \mu_E(M,E)|M\cap N|^k.
\]
Indeed, a subgroup containing $H$ can contain all the lifts only if
$MN=E$, and then each entry has $|M\cap N|$ lifts in it.  The count is
independent of the tuple, and is positive for the image of an existing
relative generating tuple.
\end{proof}

For a finite group $E$, its \emph{Frattini subgroup} $\Phi(E)$
is the intersection of all maximal proper subgroups of $E$
(with $\Phi(1)=1$).
A finite group $L$ is \emph{monolithic} if it has a unique minimal
normal subgroup, called its \emph{monolith}. We call $L$
\emph{monolithic primitive} if, in addition, $\Phi(L)=1$,
and write $A$ for its monolith.
Its crown-based power is
\begin{equation}\label{eq:crown-definition}
 L_t=\{(l_1,\ldots,l_t)\in L^t:l_1A=\cdots=l_tA\}.
\end{equation}
It has socle $A^t$ and quotient $L/A$.  If $A$ is abelian, then
$L=A\rt (L/A)$ with $A$ a faithful irreducible $\F_p[L/A]$-module; put
$q_A=|\End_{L/A}(A)|$ and $\tau_A=\dim_{\F_{q_A}}A$.
Otherwise $A$ is a power of a nonabelian simple group, and conjugation
embeds $L$ in $\Aut(A)\le\Sym_{|A|}$.

\begin{lemma}\label{lem:relative-crown}
Let $E$ be a finite group and $H\le E$. For every integer
$1\le r\le d_H(E)$, there is an epimorphism $\pi:E\onto L_t$
onto a crown-based power such that
$d_{\pi(H)}(L_t)\ge r$, while every proper quotient of $L_t$ can be
generated by the image of $\pi(H)$ and at most $r-1$ additional elements.
\end{lemma}
\begin{proof}
Choose an epimorphism $\pi:E\onto Q$ with $|Q|$ least subject to
$d_{\pi(H)}(Q)\ge r$. Every proper quotient of $Q$ has relative generator
number at most $r-1$, so it remains to show that $Q$ is a crown-based power.
Factoring out $\Phi(Q)$ does not change the relative generator number:
if $J\le Q$ and $J\Phi(Q)=Q$, then $J=Q$. Thus minimality gives $\Phi(Q)=1$.

The crown splitting lemmas \cite[Lemmas~3 and~4]{LucchiniChebotarev}
give normal subgroups $R,D\normal Q$ with $D\ne1$, with $Q/R$ a
crown-based power, and with the following property for every $J\le Q$:
\[
 JR=JD=Q\quad\Longrightarrow\quad J=Q.
\]
Suppose $R\ne1$. Since $Q/R$ is a proper quotient, choose
$x_1,\ldots,x_{r-1}\in Q$ (an empty tuple if $r=1$) such that
\[
 Q=\langle\pi(H),x_1,\ldots,x_{r-1}\rangle R.
\]
Their images generate $Q/RD$ together with the image of $\pi(H)$.
The proper quotient $Q/D$ also has relative generator number at most
$r-1$. Apply \cref{lem:relative-lifting} to $Q/D\onto Q/RD$:
there are $y_1,\ldots,y_{r-1}\in R$ such that
\[
 Q=\langle\pi(H),x_1y_1,\ldots,x_{r-1}y_{r-1}\rangle D.
\]
Set $J=\langle\pi(H),x_1y_1,\ldots,x_{r-1}y_{r-1}\rangle$.
Multiplying the $x_i$ by elements of $R$ leaves their images in $Q/R$
unchanged, so $JR=Q$; the last display gives $JD=Q$. Hence $J=Q$,
contrary to $d_{\pi(H)}(Q)\ge r$. Therefore $R=1$, and $Q$ is a
crown-based power, as required.
\end{proof}

\begin{lemma}\label{lem:relative-multiplicity}
Let $L_t$ be a crown-based power with socle $A^t$, let $H\le L_t$,
and let $r\ge4$ be an integer. Suppose that $d_H(L_t)\ge r$,
while every proper quotient of $L_t$ can be generated by the image
of $H$ and at most $r-1$ additional elements. Then
\begin{equation}\label{eq:relative-multiplicity}
 t>(r-3)\tau_A\quad(A\text{ abelian}),\qquad
 t>|A|^{r-3}\quad(A\text{ nonabelian}).
\end{equation}
In particular, $t>1$ and $d(L)\le r-1+d(H)$,
where $d(G)$ denotes $d_\Gamma(G)$ with $\Gamma$ taken to be
the trivial group.
\end{lemma}
\begin{proof}
Suppose first that $A$ is abelian. Fix $r-1$ elements of $L/A$
which generate this quotient together with the image of $H$.
In $L_t=A^t\rt(L/A)$, choose a lift $x_i$ of each fixed element,
independently and uniformly among its $|A^t|=q_A^{\tau_A t}$ lifts.
All probabilities below refer to this random tuple
$(x_1,\ldots,x_{r-1})$.

If these lifts do not generate $L_t$ together with $H$, they lie in
a maximal subgroup $M$ containing $H$. Since their images together
with the image of $H$ generate $L_t/A^t$, this subgroup satisfies
$MA^t=L_t$. Then $M\cap A^t$ is a maximal $(L/A)$-submodule of $A^t$, with
quotient $A$, and $M/(M\cap A^t)$ is a complement to $A$ in
$A\rt(L/A)$. There are $(q_A^t-1)/(q_A-1)<q_A^t$ such submodules.
For each one, the number of complements is
\[
 \frac{|A|}{|A^{L/A}|}\,|H^1(L/A,A)|
 \le |A|\,|H^1(L/A,A)|<q_A^{2\tau_A},
\]
by \cite[Theorem~A]{AschbacherGuralnick}; when $L/A=1$, the
last inequality is immediate. Thus fewer than $q_A^tq_A^{2\tau_A}$ maximal
subgroups $M$ can occur.

For each such $M$, every coset of $A^t$ meets $M$ in exactly
$|M\cap A^t|$ elements. Since $A^t/(M\cap A^t)\simeq A$,
uniformity and independence give
\[
 \Prob(x_i\in M)
 =\frac{|M\cap A^t|}{|A^t|}=\frac1{|A|}=q_A^{-\tau_A},
 \Prob(x_1,\ldots,x_{r-1}\in M)=q_A^{-\tau_A(r-1)}.
\]
Failure to generate is the union of these events over the maximal
subgroups just counted. Consequently, if $t\le(r-3)\tau_A$, the union bound gives
\[
 \Prob\bigl(\langle H,x_1,\ldots,x_{r-1}\rangle\ne L_t\bigr)
 \le\sum_M\Prob(x_1,\ldots,x_{r-1}\in M)
 <q_A^tq_A^{2\tau_A}q_A^{-\tau_A(r-1)}
 =q_A^{t-(r-3)\tau_A}\le1.
\]
The failure probability is strictly less than one, so some choice of
lifts generates $L_t$ together with $H$, contradicting $d_H(L_t)\ge r$.

Now suppose that $A$ is nonabelian, and choose a coordinate minimal
normal subgroup $N\simeq A$. The hypothesis on $L_t/N$ gives
\[
 L_t=\langle H,x_1,\ldots,x_{r-1}\rangle N.
\]
Choose generators $y_1,\ldots,y_{d(H)}$ of $H$ and append them to obtain
\[
 (x_1,\ldots,x_{r-1},y_1,\ldots,y_{d(H)}).
\]
If $t\le |A|^{r-3}$, then
\[
 h:=\left\lceil\frac85+\log_{|A|} t\right\rceil
 \le\left\lceil\frac85+r-3\right\rceil=r-1.
\]
By \cite[Proposition~6]{LucchiniThakkar}, there are
$n_1,\ldots,n_h\in N$ such that
\[
 L_t=\langle x_1n_1,\ldots,x_hn_h,
 x_{h+1},\ldots,x_{r-1},y_1,\ldots,y_{d(H)}\rangle.
\]
Only the first $h$ entries $x_1,\ldots,x_h$ have been modified.
Thus the generators of $H$ remain fixed, and $r-1$ additional
elements generate $L_t$, again a contradiction.

Since $r\ge4$, both bounds imply $t>1$. Hence $L$ is a proper
coordinate quotient of $L_t$ and can be generated by the image of
$H$ and at most $r-1$ additional elements. This gives the bound on $d(L)$.
\end{proof}

\subsubsection{Extraction of quotients with fixed multiplicity}\label{subsec:fixed-multiplicity}

We apply the relative crown lemmas to $G\rt\Gamma$ and then project
to a prescribed multiplicity. The resulting quotients can be counted
and their probabilities bounded uniformly. This projection preserves
nonabelianness, as needed for the final comparison with the abelianization.

\begin{lemma}\label{lem:threshold-witness}
Let $G$ be a finite $\Gamma$-group with $G_\Gamma=1$,
$(|G|,|\Gamma|)=1$ and $d_\Gamma(G)\ge r\ge4$.
There is an epimorphism $\pi:G\rt\Gamma\onto L_t$ onto a
crown-based power such that
\begin{equation}\label{eq:threshold-witness}
 L_t=\pi(G)\rt \pi(\Gamma),\qquad A^t\le \pi(G),\qquad d_{\pi(\Gamma)}(L_t)\ge r,
\end{equation}
and every proper quotient has relative rank at most $r-1$.
Moreover $(\pi(G))_\Gamma=1$, $L_t$ and $L/A$ are normally generated by the image
of $\pi(\Gamma)$, and
\begin{equation}\label{eq:base-rank}
 d(L/A)\le r-1+d(\Gamma),\qquad
 t>(r-3)\tau_A\ \ (A\text{ abelian}),\qquad
 t>|A|^{r-3}\ \ (A\text{ nonabelian}).
\end{equation}
In particular, $t>1$ and $d(L)\le r-1+d(\Gamma)$.
\end{lemma}
\begin{proof}
Apply \cref{lem:relative-crown} with $E=G\rt\Gamma$ and $H=\Gamma$
to obtain $\pi:G\rt\Gamma\onto L_t$. Coprimality gives
$\pi(G)\cap\pi(\Gamma)=1$, so $L_t=\pi(G)\rt\pi(\Gamma)$.
Applying \cref{lem:relative-multiplicity} with $H=\pi(\Gamma)$ gives
the multiplicity bounds, $t>1$, and
$d(L/A)\le d(L)\le r-1+d(\Gamma)$.

Set $D=A^t\cap\pi(G)$. The semisimple structure of $A^t$ gives
$D\simeq A^j$ and $A^t/D\simeq A^{t-j}$ for some $0\le j\le t$.
Since the latter embeds in $L_t/\pi(G)$, coprimality forces
$j=t$ or $j=0$. If $j=0$, then $D=1$, so the normal subgroups
$\pi(G),A^t$ centralize each other; $C_{L_t}(A^t)\le A^t$ then gives
$\pi(G)=1$, impossible. Hence $j=t$ and $A^t\le\pi(G)$.
The relation $(\pi(G))_\Gamma=1$ descends from $G$ and is equivalent
to normal generation of $L_t$ by $\pi(\Gamma)$.  In the abelian case $L/A\ne1$,
since otherwise $\pi(\Gamma)=1$ and $(\pi(G))_\Gamma=\pi(G)\ne1$.
\end{proof}

\begin{definition}\label{def:extracted-witness}
Fix $r\ge4$. A \emph{witness} is a tuple $(G,L,A,\pi,r)$ satisfying
the hypotheses and conclusions of \cref{lem:threshold-witness}; in
particular, $\pi:G\rt\Gamma\onto L_t$.

To \emph{extract} a witness, set
\begin{equation}\label{eq:fixed-multiplicity}
 j=\begin{cases}
    (r-3)\tau_A, & A\text{ abelian},\\
    |A|^{r-3}, & A\text{ nonabelian}.
   \end{cases}
\end{equation}
By \eqref{eq:base-rank}, $1\le j<t$. Compose $\pi$ with the
projection $L_t\onto L_j$ onto the first $j$ coordinates, obtaining
$\pi_j:G\rt\Gamma\onto L_j$. The tuple
$(G,L,A,\pi_j,r)$ is its \emph{extracted witness}.

A witness or extracted witness $(G,L,A,\varphi,r)$ is
\emph{nonabelian} if $\varphi(G)$ is nonabelian. Its
\emph{socle type} is the type of $A$ (abelian or nonabelian).

The \emph{witness quotient} of an extracted witness is the
$\Gamma$-group $W=\pi_j(G)$, with
$\gamma\cdot w=\pi_j(\gamma)w\pi_j(\gamma)^{-1}$.
An extracted witness is a \emph{normalized witness} if
$\pi_j(\Gamma)$ is diagonal in $L_j$ and, when $A$ is abelian,
lies in the standard complement $L/A$.

For counting witness quotients, associate to each normalized witness
its \emph{normalized data}:
\begin{itemize}
\item If $A$ is abelian, the subgroup
$L/A\le\GL_{\tau_A}(q_A)$ and the homomorphism
$\Gamma\to L/A$ induced by $\pi_j(\Gamma)$.
\item If $A$ is nonabelian, the subgroup
$L\le\Sym_{|A|}$ (acting by conjugation on $A$) and the
homomorphism $\Gamma\to L$ given by a coordinate of the diagonal
image $\pi_j(\Gamma)$.
\end{itemize}
For fixed $r$, these data determine $\pi_j(G)$ up to $\Gamma$-isomorphism,
as proved in \cref{lem:normalized-data}; they do not record the source
$G$ or the map $\pi_j$.
\end{definition}

\begin{lemma}\label{lem:normalized-data}
Every extracted witness $(G,L,A,\pi_j,r)$ can be normalized by
conjugating $\pi_j$ with an element of $A^j$. If the original witness is
nonabelian, its extracted witness is nonabelian.

For fixed $q_A,\tau_A$ in the abelian case, or fixed $|A|$ in the
nonabelian case, the number of possible normalized data is at most
\begin{equation}\label{eq:base-counts}
 \begin{cases}
 q_A^{(r-1+2d(\Gamma))\tau_A^2}, & A\text{ abelian},\\
 (|A|!)^{r-1+2d(\Gamma)}, & A\text{ nonabelian}.
 \end{cases}
\end{equation}
For every normalized witness, moreover,
\begin{equation}\label{eq:automorphism-estimates}
 |\GammaAut(W)|\ge
 \begin{cases}
 |\GL_j(q_A)|\ge\rho q_A^{j^2}, & A\text{ abelian},\\
 j!, & A\text{ nonabelian},
 \end{cases}
 \qquad\rho=\eta(2).
\end{equation}
\end{lemma}
\begin{proof}
\noindent\textbf{1. Normalizing the image of $\Gamma$.}

Choose an original witness $(G,L,A,\pi,r)$ giving the extracted
witness, and let $\pi_1:G\rt\Gamma\onto L$ be the first coordinate
of $\pi$.
The map $\pi(\Gamma)\to L/A$ is injective, since
$A^t\le\pi(G)$ and $\pi(G)\cap\pi(\Gamma)=1$. Put
\[
 E:=A\pi_1(\Gamma),\qquad
 E/A=A\pi_1(\Gamma)/A\simeq\pi(\Gamma).
\]
Thus $(|A|,|E/A|)=(|A|,|\pi(\Gamma)|)=1$.
The coordinate images of $\pi(\Gamma)$ in $L$ are complements to $A$
in $E$.

\smallskip\noindent\emph{Abelian case.}
Choose the common complement to $A$ in $E$ inside the standard
complement $L/A$.

\smallskip\noindent\emph{Nonabelian case.}
Choose $\pi_1(\Gamma)$ as the common complement to $A$ in $E$.

\smallskip\noindent
In both cases, by the full conjugacy form of Schur--Zassenhaus
\cite[Theorem~6.2.1 and p.~127]{KurzweilStellmacher}, we can conjugate
each coordinate image by an element of $A$ onto the chosen complement.
For each $\gamma\in\Gamma$, the conjugated coordinates lie in the same complement
and project to the same element of $E/A$. The projection is injective
on that complement, so the coordinates are equal. Conjugating by the
tuple of the chosen elements of $A$, viewed as an element of $A^t$,
therefore makes $\pi(\Gamma)$ diagonal. Its projection to $A^j$
conjugates $\pi_j$ to a normalized extracted witness.

\medskip\noindent\textbf{2. Properties of $\pi_j(G)$.}

The subgroup $\pi_j(G)$ is a $\Gamma$-quotient of $G$, with
$(\pi_j(G))_\Gamma=1$.
It is uniquely determined as a subgroup of $L_j$: it is the unique
normal subgroup whose order is prime to $|\Gamma|$ and whose index
has all its prime divisors among those of $|\Gamma|$.
It is therefore preserved by every automorphism of $L_j$.

\smallskip\noindent\emph{Abelian case.}
Write $N=\pi_j(G)/A^j\simeq\pi(G)/A^t\le L/A$.
Then $\pi_j(G)=A^j\rt N$. If the original witness is nonabelian,
then $N\ne1$, and its faithful action on $A$ makes $\pi_j(G)$
nonabelian. Conjugation preserves this property.

\smallskip\noindent\emph{Nonabelian case.}
The subgroup $A^j\le\pi_j(G)$ already makes $\pi_j(G)$ nonabelian.

\medskip\noindent\textbf{3. Counting the possible normalized data.}

\smallskip\noindent\emph{Abelian case.}
The data are $L/A\le\GL_{\tau_A}(q_A)$ and a map $\Gamma\to L/A$.
By \eqref{eq:base-rank}, we can specify the subgroup by a generating
$(r-1+d(\Gamma))$-tuple, and the map by the images of $d(\Gamma)$
generators. Thus the number of choices is at most
\[
 |\GL_{\tau_A}(q_A)|^{r-1+2d(\Gamma)}
 \le q_A^{(r-1+2d(\Gamma))\tau_A^2}.
\]
These data determine $L_j$ and its $\Gamma$-action;
arbitrary tuples only overcount.

\smallskip\noindent\emph{Nonabelian case.}
Embed $L\le\Sym_{|A|}$ by conjugation on $A$.
The bound $d(L)\le r-1+d(\Gamma)$ lets us count generating tuples
for $L$ and the images of $d(\Gamma)$ generators under the normalized
map $\Gamma\to L$. There are at most
\[
 (|A|!)^{r-1+d(\Gamma)}(|A|!)^{d(\Gamma)}
 =(|A|!)^{r-1+2d(\Gamma)}
\]
choices. The unique minimal normal subgroup of $L$ and the prescribed
$j$ determine $L_j$ and its $\Gamma$-action. In either case, part~2
identifies $W$ as a characteristic subgroup of $L_j$, so the normalized
data determine $W$ up to $\Gamma$-isomorphism, though they do not
specify $G$ or $\pi_j$. This proves \eqref{eq:base-counts}.

\medskip\noindent\textbf{4. Bounding $|\GammaAut(\pi_j(G))|$ from below.}

\smallskip\noindent\emph{Abelian case.}
Every element of $\Aut_{L/A}(A^j)\simeq\GL_j(q_A)$ extends to an
automorphism of $L_j=A^j\rt(L/A)$ acting identically on $L/A$.
It fixes $\pi_j(\Gamma)$ pointwise and preserves $\pi_j(G)$ by part~2,
so its restriction is $\Gamma$-equivariant. Restriction is faithful
because $A^j\le\pi_j(G)$, giving
\[
 |\GammaAut(\pi_j(G))|\ge|\GL_j(q_A)|\ge\rho q_A^{j^2}.
\]

\smallskip\noindent\emph{Nonabelian case.}
Coordinate permutations act on $L_j$ and fix the diagonal
$\pi_j(\Gamma)$ pointwise. They preserve $\pi_j(G)$ by part~2
and act faithfully on $A^j\le\pi_j(G)$.
Thus $|\GammaAut(\pi_j(G))|\ge j!$, proving
\eqref{eq:automorphism-estimates}.
\end{proof}

\subsection{The quadratic-exponent gap for nonabelian witnesses with abelian socle}\label{sec:rigidity}

When $A$ is nonabelian, the large multiplicity and factorial
automorphism bound from \cref{lem:normalized-data} give the required
probability estimate in \cref{sec:estimates}. The remaining case is
a nonabelian normalized witness with abelian socle $A$, so its
quotient $W=\pi_j(G)$ is nonabelian.
To show that the probability of such a quotient is negligible,
we need its quadratic decay coefficient $\kappa(A)$, defined below, to be
strictly larger than the leading coefficient $\kappa$.

For a normalized witness with abelian socle, set
\begin{equation}\label{eq:abelian-socle-quadratic-coefficient}
 a(A)=\begin{cases}
 \tfrac12,&p\mid n\text{ and }A\simeq A^\vee
                         \text{ as }\F_p[L/A]\text{-modules},\\
 1,&\text{otherwise},
 \end{cases}
 \qquad \kappa(A)=a(A)\tau_A^2\log q_A.
\end{equation}
All duals in this subsection are over the prime field.

\begin{lemma}\label{lem:quadratic-coefficient-comparison}
We have $\kappa(A)\ge\kappa$.  Equality implies
\begin{equation}\label{eq:equality-restriction}
 A|_\Gamma\simeq V\in\cM,\qquad q_A=q_V,\quad \tau_A=\tau_V,\quad a(A)=a_V.
\end{equation}
\end{lemma}
\begin{proof}
The restriction has a nontrivial constituent: otherwise $\pi_j(\Gamma)$ acts
trivially, and its normal generation makes the faithful group $L/A\ne1$
trivial.  Let $\mathcal R$ be the image of the algebra homomorphism
$\F_p[\Gamma]\to\End_{\F_{q_A}}(A)$ induced by the $\Gamma$-action
on $A$. For an $\F_p[\Gamma]$-module $W$, let
$\operatorname{Supp}(W)$ denote the set of isomorphism classes of its
irreducible constituents, counted without multiplicities. Then
\[
 \mathcal R\simeq\prod_{V\in\operatorname{Supp}(A|_\Gamma)}
 M_{\tau_V}(q_V)\ \subseteq M_{\tau_A}(q_A),\qquad
 \sum_{V\in\operatorname{Supp}(A|_\Gamma)}\tau_V^2\log q_V
      \le \tau_A^2\log q_A.
\]
For the trivial module, if present, $q_V,\tau_V$ have the same
endomorphism-field meanings.  If $a(A)=1$, any nontrivial constituent gives the bound.
If $a(A)=1/2$, use a nontrivial self-dual constituent when one exists.
Otherwise self-duality gives a distinct dual pair $V,V^\vee$; the image
algebra has at least two factors, hence $\mathcal R\subsetneq M_{\tau_A}(q_A)$, and
\[
 \kappa(A)=\tfrac12\log|M_{\tau_A}(q_A)|>
 \tfrac12\log|\mathcal R|\ge\kappa_V\ge\kappa.
\]
Equality in the other cases forces $\mathcal R=M_{\tau_A}(q_A)$ and
$a(A)=a_V$.  The natural module is then irreducible over
$\F_p[\Gamma]$ with endomorphism field $\F_{q_A}$, proving the assertion.
\end{proof}

\begin{proposition}\label{prop:quadratic-coefficient-rigidity}
For a normalized witness with abelian socle, if $\kappa(A)=\kappa$,
then $\pi_j(G)=A^j$.
Consequently there is $\kappa_+>\kappa$, depending only on $\Gamma,n$,
such that every nonabelian normalized witness with abelian socle has
$\kappa(A)\ge\kappa_+$.
\end{proposition}
\begin{proof}
Put $Q=\pi_j(G)$ and $H=\pi_j(\Gamma)\le L/A$.
We show that $N:=Q/A^j\le L/A$ is trivial, which gives $Q=A^j$.
We separate the cases where $|\Gamma|$ is even and odd.
An irreducible representation is \emph{primitive} if its space
admits no direct-sum decomposition into two or more nonzero subspaces
permuted by the group. In the even case, we first show that
$\bigl(A\otimes_{\F_{q_A}}\overline{\F}_p\bigr)|_\Gamma$
is a primitive representation of $\Gamma$.
Since $N$ has odd order and is solvable, Clifford theory and
the Fong--Swan theorem then force $N$ to act by scalars.
In the odd case, $H$ is solvable, and a cyclic quotient produces
a representation $W$ with $\kappa_W<\kappa$ unless $\tau_A=1$.
In both cases, normal generation by $H$ then forces $N=1$.

In both cases, we will use the divisibility $\tau_A\mid|\Gamma|$.
By \cref{lem:quadratic-coefficient-comparison}, $A|_\Gamma\simeq V\in\cM$, and the
restriction of $A$ to $H$ is absolutely irreducible over
$\F_{q_A}$: its image algebra is $M_{\tau_A}(q_A)$. Set
$\overline A=A\otimes_{\F_{q_A}}\overline{\F}_p$.
Nonmodular character theory gives $\tau_A=\tau_V\mid|\Gamma|$;
see \cite{IsaacsBook}. Recall also that
\[
 L/A=N\rt H\le\GL_{\tau_A}(q_A),\qquad (|N|,|\Gamma|)=1,
\]
and that $H$ normally generates $L/A$.

\medskip\noindent\textbf{The case $|\Gamma|$ even.}

\smallskip\noindent\emph{Equality $\kappa(A)=\kappa$ forces primitivity.}
Suppose that $\overline A=B_1\oplus\cdots\oplus B_k$, with $k>1$,
is a decomposition into nonzero subspaces permuted by $\Gamma$.
Irreducibility makes this permutation action transitive, so
$2\le k\le \tau_A$ and $k\mid|\Gamma|$.
Let $\Gamma$ act on $\F_p^k$ by the permutation of the $B_i$:
if $\gamma(B_i)=B_j$, then $\gamma(e_i)=e_j$ for the standard basis vectors.
This action preserves the augmentation module
\[
 M=\{(x_1,\ldots,x_k)\in\F_p^k:x_1+\cdots+x_k=0\}.
\]
The $\Gamma$-module $M$ is semisimple since $p\nmid|\Gamma|$.
Transitivity makes every $\Gamma$-fixed vector of $\F_p^k$ constant,
and such a vector belongs to $M$ only if it is zero, since $p\nmid k$.
Thus $M$ has no trivial constituent.

The dot product is $\Gamma$-invariant because this action permutes coordinates.
Its restriction to $M$ is nondegenerate: a vector of $M$ orthogonal to
all of $M$ must have constant coordinates, and hence be zero since $p\nmid k$.
Consequently the map
\[
 M\longrightarrow M^\vee,\qquad
 x\longmapsto\left(y\longmapsto\sum_{i=1}^k x_i y_i\right)
\]
is a $\Gamma$-equivariant isomorphism, so $M$ is self-dual as a
$\Gamma$-representation.

For a constituent $W$ with $\dim_{\F_p}W=d$, writing $q_W=p^f$ gives
$d=f \tau_W$ and hence $\tau_W^2\log q_W\le d^2\log p$.
If $a(A)=1$, choose any constituent of $M$.
If $a(A)=1/2$, then $p\mid n$: choose a self-dual constituent when
one exists, so that $a_W=1/2$; otherwise choose a member of a
distinct dual pair, of dimension at most $(k-1)/2$.
In every case,
\[
 \kappa\le\kappa_W\le a(A)(k-1)^2\log p
        <a(A)\tau_A^2\log q_A=\kappa,
\]
a contradiction. Thus $\overline A$ admits no such decomposition.

\smallskip\noindent\emph{Applying Clifford theory and Fong--Swan.}
The group $N$ has odd order, so it is solvable by \cite{FeitThompson}.
By Clifford theory, $\overline A|_N$ is completely reducible, and its
isotypic components are permuted by $H$.
The primitivity just proved forces a single isotypic component, so
\[
 \overline A|_N\simeq V_0^{\oplus e}
\]
for an irreducible $\overline{\F}_p[N]$-module $V_0$.
Put $d=\dim_{\overline{\F}_p}V_0$; then $\tau_A=e d$.
Since $N$ is solvable, the Fong--Swan theorem
\cite[Theorem~10.1]{NavarroBook} lifts the Brauer character of $V_0$
to an ordinary irreducible character of the same degree.
Ordinary irreducible character degrees divide the group order, so
$d\mid|N|$. But also $d\mid \tau_A\mid|\Gamma|$, and coprimality gives $d=1$.
All summands are copies of the same one-dimensional representation,
so every element of $N$ acts by a scalar on $\overline A$.
Faithfulness makes $N$ central in $L/A$, and hence
$L/A=N\times H$.
Thus $H$ is normal; since it normally generates $L/A$,
we obtain $L/A=H$ and $N=1$.

\medskip\noindent\textbf{The case $|\Gamma|$ odd.}
We claim that $\tau_A=1$. Otherwise $\tau_A\mid|\Gamma|$ gives $\tau_A\ge3$.
The nontrivial image $H$ in $\GL_{\tau_A}(q_A)$ is solvable
by \cite{FeitThompson}, so it has a quotient $C_\ell$ for a prime $\ell\ne p$.
Since $\ell\mid|\GL_{\tau_A}(q_A)|$, we have $\ell\mid q_A^i-1$
for some $i\le \tau_A$.
Pulling back a nontrivial simple $\F_p[C_\ell]$-module gives
$W\in\cI$ with $\tau_W=1$ and $q_W\le q_A^i\le q_A^{\tau_A}$.
Consequently
\[
 \kappa\le\kappa_W\le \tau_A\log q_A
 <\tfrac12\tau_A^2\log q_A\le\kappa(A)=\kappa,
\]
a contradiction. Thus $\tau_A=1$, and $L/A\le\GL_1(q_A)$ is abelian.
Its normal generation by $H$ again gives
$L/A=H$ and $N=1$.

\medskip\noindent\textbf{The strict inequality $\kappa(A)>\kappa$.}
We have proved $\pi_j(G)=A^j$ whenever $\kappa(A)=\kappa$.
By \cref{lem:quadratic-coefficient-comparison}, every nonabelian normalized witness
with abelian socle therefore has $\kappa(A)>\kappa$.
The numerical set
\[
 \{a \tau_A^2\log q_A:q_A\text{ a prime power},\ \tau_A\ge1,\ a\in\{1/2,1\}\}
\]
is locally finite.  Its least member greater than $\kappa$ can be
chosen as $\kappa_+$.
\end{proof}

\subsection{Bounding the difference between the full and abelian \texorpdfstring{\mbox{generator rank tails}}{generator rank tails}}\label{sec:estimates}

We now bound the difference between $\Prob(d_\Gamma(X)\ge r)$ and
$\Prob(d_\Gamma(X^{\ab})\ge r)$. The goal is to obtain, for some
$\delta>0$ depending only on $\Gamma,n$,
\[
 0\le\Prob(d_\Gamma(X)\ge r)-\Prob(d_\Gamma(X^{\ab})\ge r)
 \ll_{\Gamma,n} e^{-(\kappa+\delta)r^2}.
\]
If $d_\Gamma(X)\ge r>d_\Gamma(X^{\ab})$, the crown reduction and
projection in \cref{sec:relative} produce a nonabelian normalized
witness with source a finite $\Gamma$-quotient of $X$.
It therefore suffices to bound the probability of such a witness.
Every such witness provides a nonabelian witness quotient of $X$.

We convert the structural conclusions into probability estimates
using the quotient moments, and sum over normalized data parametrizing
possible nonabelian witness quotients. For abelian socles, the inequality $\kappa(A)\ge\kappa_+>\kappa$ from
\cref{sec:rigidity} makes the sum smaller than the leading elementary
abelian tail. For nonabelian socles, the exponential multiplicity and
the factorial automorphism bound give a still smaller contribution.

Constants below may depend on $\Gamma,n$ and be enlarged.
We first give a second-cohomology bound used for the base of either
kind of crown.

\begin{lemma}\label{lem:base-H2}
For every finite abelian group $T$, there is a constant $\theta_T$ such that,
for every finite group $J$ and every $J$-action on $T$,
\begin{equation}
 |H^2(J,T)|\le |J|^{\theta_T(d(J)+1)}.                    \label{eq:base-cohomology}
\end{equation}
\end{lemma}

\begin{proof}
The case $J=1$ is immediate. Otherwise, set
$C_{\mathrm{abs}}=24/\log 2$.
By \cite[Theorem~13.9]{GuralnickKantorKassabovLubotzky} and the
chief-factor bounds stated immediately after it, there is a free
profinite presentation $1\to R\to F\to J\to1$ in which $R$ is normally
generated by at most
\[
 d(J)+20+(d(J)+2)\log_2|J|
 \le C_{\mathrm{abs}}d(J)\log|J|
\]
elements. The inequality uses $d(J)\ge1$ and $\log_2|J|\ge1$.
Inflate the action on $T$ to $F$. Since $R$ acts trivially and
$H^2(F,T)=0$, the five-term sequence makes $H^2(J,T)$ a quotient of
$\Hom_{\mathrm{cont}}(R,T)^F$. Such an equivariant homomorphism is
determined by its values on the chosen normal generators. Hence
\[
 |H^2(J,T)|\le |T|^{C_{\mathrm{abs}} d(J)\log|J|}
              =|J|^{C_{\mathrm{abs}} d(J)\log|T|}.
\]
Thus we may take $\theta_T=C_{\mathrm{abs}}\log|T|$.
\end{proof}

For $T=\Z/n$, write $\theta_n=\theta_{\Z/n}$.

\begin{lemma}\label{lem:witness-H2}
For a normalized witness, with $j$ as in
\eqref{eq:fixed-multiplicity}, inflation gives
\[
 H^2(\pi_j(G)\rt\Gamma,\Z/n)\simeq H^2(L_j,\Z/n).
\]
If $A$ is abelian, then
\begin{equation}\label{eq:abelian-witness-H2}
 |H^2(L_j,\Z/n)|
 \le q_A^{(1-a(A))(j^2+j)+\tau_A j+\theta_n(r+d(\Gamma))\tau_A^2}.
\end{equation}
If $A$ is nonabelian, then
\begin{equation}\label{eq:nonabelian-witness-H2}
 |H^2(L_j,\Z/n)|
 \le (|A|!)^{\theta_n(r+d(\Gamma))}|A|^{3\theta_nj}.
\end{equation}
\end{lemma}
\begin{proof}
The kernel of $\pi_j(G)\rt\Gamma\to L_j=\pi_j(G)\rt\pi_j(\Gamma)$ has order dividing
$|\Gamma|$, hence prime to $n$. Inflation identifies their cohomology
with coefficients $\Z/n$.

\emph{Abelian socle.}
Let $A$ have characteristic $p$. Coefficients prime to $p$
only contribute base cohomology. For $p^e\Vert n$, the degree-two terms
of the Lyndon--Hochschild--Serre spectral sequence give
\[
 |H^2(L_j,\Z/p^e)|\le |H^2(L/A,\Z/p^e)|\,
 |H^1(L/A,(A^\vee)^j)|\,|H^2(A^j,\Z/p^e)^{L/A}|.
\]
The middle factor is at most $q_A^{\tau_A j}$ by \cite[Theorem~A]{AschbacherGuralnick}.
The universal coefficient sequence has no $(L/A)$-invariants in its
$\Ext^1$ term, since $A$ is nontrivial irreducible module. It therefore gives
an injection of the last group into
$\Hom(\bigwedge^2_{\F_p}A^j,\F_p)^{L/A}$.
Schur's lemma and the exterior-power decomposition bound
this by $q_A^{\binom j2+j}$ if $A$ is self-dual, and by $1$ otherwise.
Applying \eqref{eq:base-cohomology}, with $|L/A|\le q_A^{\tau_A^2}$ and
$d(L/A)+1\le r+d(\Gamma)$, gives \eqref{eq:abelian-witness-H2}.

\emph{Nonabelian socle.}
Write $A=S^a$ with $S$ simple. Since $A^j$ is perfect, the same
spectral sequence has no middle $H^1$-term and gives
\[
 |H^2(L_j,\Z/n)|\le |H^2(L/A,\Z/n)|\,
                         |H^2(A^j,\Z/n)^{L/A}|.
\]
Every nonabelian finite simple group is generated by two elements;
see \cite[p.~471]{GuralnickKantorKassabovLubotzky}. Thus
\eqref{eq:base-cohomology}, universal coefficients, and the product
formula for the Schur multiplier of perfect groups give
\[
 |H^2(A^j,\Z/n)^{L/A}|\le |H^2(S,\Z/n)|^{aj}
                       \le |A|^{3\theta_nj}.
\]
Applying \eqref{eq:base-cohomology} to $L/A$, with $|L/A|\le |A|!$ and
$d(L/A)+1\le r+d(\Gamma)$, gives \eqref{eq:nonabelian-witness-H2}.
\end{proof}

\begin{lemma}\label{lem:witness-count}
Set
\[
 C_0=3+(d(\Gamma)+3)\theta_n
    =3+\frac{24}{\log 2}(d(\Gamma)+3)\log n.
\]
For a normalized witness, with $j$ as in
\eqref{eq:fixed-multiplicity}, the expected number of open normal
$\Gamma$-stable subgroups of $X$ with quotient $\Gamma$-isomorphic to
$\pi_j(G)$ satisfies
\[
 \begin{aligned}
 \Exp\#\{N\normal X:N\text{ open and }\Gamma\text{-stable},
                    \,X/N\simeq_\Gamma\pi_j(G)\}\le\begin{cases}
 \rho^{-1}q_A^{-a(A)j^2+C_0(\tau_A j+r \tau_A^2)},&A\text{ abelian},\\
 (|A|!)^{C_0r}\,|A|^{C_0j}/j!,&A\text{ nonabelian}.
 \end{cases}
 \end{aligned}
\]
\end{lemma}
\begin{proof}
By \eqref{eq:underlying-moment}, \cref{lem:witness-H2}, and
$w_U(\pi_j(G))\ge1$,
\[
 \frac{\Exp|\GammaSur(X,\pi_j(G))|}{|\GammaAut(\pi_j(G))|}
 \le \frac{|H^2(L_j,\Z/n)|\,}
              {|\GammaAut(\pi_j(G))|}.
\]

\emph{Abelian socle.}
Use \eqref{eq:abelian-witness-H2} and
\eqref{eq:automorphism-estimates}. Since $\tau_A\ge1$ and
$1-a(A)\le1/2$, we obtain
\[
 \frac{\Exp|\GammaSur(X,\pi_j(G))|}{|\GammaAut(\pi_j(G))|}
 \le\rho^{-1}q_A^{-a(A)j^2+3\tau_A j+
                  [1+\theta_n(r+d(\Gamma))]\tau_A^2}.
\]
The inequalities $C_0\ge3$ and
$1+\theta_n(r+d(\Gamma))\le C_0r$ prove the abelian case.

\emph{Nonabelian socle.}
Use \eqref{eq:nonabelian-witness-H2} and
$|\GammaAut(\pi_j(G))|\ge j!$ to obtain
\[
 \frac{\Exp|\GammaSur(X,\pi_j(G))|}{|\GammaAut(\pi_j(G))|}
 \le (|A|!)^{1+\theta_n(r+d(\Gamma))}
       \frac{|A|^{(1+3\theta_n)j}}{j!}.
\]
Since $1+\theta_n(r+d(\Gamma))\le C_0r$ and $1+3\theta_n\le C_0$,
this proves the nonabelian case.
\end{proof}

\begin{proposition}\label{prop:negligible-witnesses}
There exists $\delta>0$, depending only on $\Gamma,n$, such that
for every $r\ge4$, the probability that some finite $\Gamma$-quotient $G$ of $X$
is the source of a nonabelian normalized witness $(G,L,A,\pi_j,r)$
is $O_{\Gamma,n}(e^{-(\kappa+\delta)r^2})$.
\end{proposition}
\begin{proof}
For each choice of normalized data, let $W$ be the associated
witness quotient. The probability that $X$ has $W$ as a quotient is
at most the expected number of its quotient kernels. Having a
normalized witness implies having the associated $W$ as a quotient,
so a union bound over normalized data suffices, even if several
witnesses give the same $\Gamma$-isomorphism type of $W$.

Set
\[
 C_1=8+4C_0+d(\Gamma),\qquad C_2=C_0+1+2d(\Gamma),
\]
where $C_0$ is the explicit constant in \cref{lem:witness-count}.

In the abelian-socle case, multiply the inequality in \cref{lem:witness-count} by
the count \eqref{eq:base-counts} and substitute $j=(r-3)\tau_A$.
For fixed $q_A,\tau_A,a$ the result is at most
\[
 \rho^{-1}q_A^{\tau_A^2[-a(r-3)^2+C_0(2r-3)+r-1+2d(\Gamma)]}
 \le
 \rho^{-1}\exp\{-a\tau_A^2\log q_A
                   [r^2-C_1r]\}.
\]
Here we used $r\ge4$ and $a\ge1/2$: these give
$2d(\Gamma)\le d(\Gamma)r/2$ and
$2C_0+1+d(\Gamma)/2\le a(4C_0+2+d(\Gamma))$.
For witnesses with nonabelian witness quotient $W$ and abelian monolith $A$, $\kappa(A)\ge\kappa_+>\kappa$.
For $r^2-C_1r>4$, summation is bounded by
\[
 \rho^{-1}e^{-\kappa_+[r^2-C_1r-4]}
 {}\times
 \sum_{q_A\ge2}\sum_{\tau_A\ge1}\sum_{a\in\{1/2,1\}}q_A^{-4a\tau_A^2}
 \ll e^{-\kappa_+r^2+O(r)}.
\]
Here the sum is convergent even over all integers $q_A\ge2$, since
$4a\tau_A^2\ge2\tau_A$. This has the required bound for any fixed
$0<\delta<\kappa_+-\kappa$.

For nonabelian socles, count directly by $|A|$. The estimate
\eqref{eq:base-counts} and the inequality in \cref{lem:witness-count}
bound the total for each $|A|$ by
\[
 (|A|!)^{C_2r}
 \frac{|A|^{C_2j}}{j!},
 \qquad j=|A|^{r-3}.
\]
Indeed, the factorial exponent from these bounds is
$(C_0+1)r-1+2d(\Gamma)\le C_2r$,
and the exponent of $|A|$ is at most $C_2j$.
For all sufficiently large $r$, uniformly for $|A|\ge2$, Stirling's
inequality gives
\[
 \begin{aligned}
 &\log\left((|A|!)^{C_2r}
            \frac{|A|^{C_2j}}{j!}\right)
 \le C_2r|A|\log|A|{}-j\bigl((r-3-C_2)\log|A|-1\bigr)\le -j\log|A|.
 \end{aligned}
\]
For the last inequality, we require $r$ to be sufficiently large that
\[
 C_2r|A|
 \le(r-C_2-6)|A|^{r-3},
\]
uniformly in $|A|\ge2$.
Consequently, summing over all integer values of $|A|\ge2$ gives
\[
 \sum_{|A|\ge2}|A|^{-|A|^{r-3}}
 \le 2\,2^{-2^{r-3}}=O(e^{-r^3}),
\]
using $|A|^{r-3}\ge2^{r-3}+|A|-2$. This decays faster than $e^{-cr^2}$
for every fixed $c>0$. All unions are countable, and the bounds are
uniform over $U$ satisfying the fixed-vector condition stated in the
introduction, since only
$w_U(\pi_j(G))^{-1}\le1$ was used.
\end{proof}

\subsection{Completion of the sharp-tail proof}\label{sec:conclusion}

We combine the exact abelian calculation with the negligible-witness
estimate. If the full group needs at least $r$ generators but its
abelianization does not, the crown reduction produces a nonabelian
normalized witness; its existence is bounded via the witness quotient. This transfers
the abelian asymptotic to the full group. We then select the least
linear exponent among the representations in $\cM$ to obtain $C$
and $\zeta$, and deduce the point-mass asymptotic.

\begin{proof}[Proof of \cref{thm:sharp-tail}]
Work on the full-measure set of \cref{prop:underlying-moment}.

\noindent\textbf{Bounding the difference between $\Prob(d_\Gamma(X)\ge r)$ and
$\Prob(d_\Gamma(X^{\ab})\ge r)$.}
If $d_\Gamma(X)\ge r>d_\Gamma(X^{\ab})$, compactness gives a finite
quotient $G$ with $d_\Gamma(G)\ge r$. Apply
\cref{lem:threshold-witness} to obtain a witness
$(G,L,A,\pi,r)$. Its image $\pi(G)$ is nonabelian, since otherwise
the epimorphism $X\onto_\Gamma\pi(G)$ would factor through
$X^{\ab}$. By \cref{lem:normalized-data}, the witness has a
nonabelian normalized extraction. Applying
\cref{prop:negligible-witnesses} gives
\begin{equation}\label{eq:comparison}
 0\le\Prob(d_\Gamma(X)\ge r)-\Prob(d_\Gamma(X^{\ab})\ge r)
       \ll e^{-(\kappa+\delta)r^2}.
\end{equation}

Together with \eqref{eq:comparison}, \cref{prop:abelian-tail} gives
\eqref{eq:finite-sum-answer}: the finite sum is at least
$e^{-\kappa r^2-O(r)}$, so the comparison error is absorbed into its
relative error $O(e^{-\eps r})$.
After decreasing $\eps$ so that $\eps\le2\kappa$ and
$\eps\le\zeta_V-\zeta$ for every $V\in\cM$ with $\zeta_V>\zeta$,
expanding $k_V(r)=\tau_Vr-(\tau_V-1)$ and using the definitions of
$C$ and $\zeta$ identifies this sum as
$Ce^{-\kappa r^2-\zeta r}(1+O(e^{-\eps r}))$: the terms with
$\zeta_V>\zeta$ are exponentially smaller and the others add to $C$.
This proves \eqref{eq:main-answer}.
The ratio of the tail at $r+1$ to that at $r$ is $O(e^{-2\kappa r})$,
so subtraction gives the same asymptotic for the point mass.
\end{proof}

\section{Almost-sure finite presentation and deficiency}
\label{sec:def-presentation}

We retain the notation of the introduction, with $(n,|\Gamma|)=1$
and $U$ satisfying the fixed-vector condition stated there.
In this section we prove \cref{thm:bounded-deficiency}. All cohomology is
continuous.

The additional input is Sawin and Wood's inequalities for finite-level
quotients that occur with positive probability. They bound the difference
between second and first cohomology;
Liu's relation formula then converts this difference into the estimate
\eqref{eq:def-bound}.

\subsection{Finite quotients of positive probability}\label{subsec:positive-quotients}

A \emph{level} is a set of isomorphism classes of finite $\Gamma$-groups
generated by finitely many groups under quotients and finite subdirect
products \cite[Section~1.1]{SawinWood}. A subdirect product is a
$\Gamma$-stable subgroup of a product whose projection onto each factor
is surjective.

Write $\mathcal P_{\Gamma,n}$ for the sample space of small profinite
$n$-oriented $\Gamma$-groups, and write
$\boldsymbol X=(X,s_X)$ for an oriented sample. For a level
$\mathcal L$, the maximal $\mathcal L$-quotient of $X$ is the quotient in
$\mathcal L$ through which every continuous $\Gamma$-equivariant
quotient map from $X$ to a group in $\mathcal L$ factors. It exists
and is finite for small $X$; see \cite[Section~1.1]{SawinWood}.
With its induced orientation, denote it by $\boldsymbol X^{\mathcal L}$. For a finite
oriented $\Gamma$-group $\boldsymbol H=(H,s_H)$, put
\[
 C(\mathcal L,\boldsymbol H)
   :=\{\boldsymbol X:
          \boldsymbol X^{\mathcal L}\simeq\boldsymbol H\}.
\]
We use the completion of $\nu_{\Gamma,n,U}$.

\begin{lemma}\label{lem:def-positive}
There is a measurable set $\Omega_+\subseteq\mathcal P_{\Gamma,n}$
of measure one such that every $C(\mathcal L,\boldsymbol H)$ containing a point of
$\Omega_+$ has positive measure. If $\boldsymbol X\in\Omega_+$ and
$\boldsymbol H=\boldsymbol X^{\mathcal L}$, then $H_\Gamma=1$.
\end{lemma}

\begin{proof}
There are countably many finite $\Gamma$-groups up to isomorphism,
and hence countably many levels. Each finite group has only finitely
many $n$-orientations. Thus the collection of $C(\mathcal L,\boldsymbol H)$ is countable.
Remove the union of its zero-measure members to obtain $\Omega_+$.
The final assertion follows from
\cite[Theorem~1.8]{SawinWood}, which assigns measure zero to a
$C(\mathcal L,\boldsymbol H)$ whose quotient has $H_\Gamma\ne1$.
\end{proof}

Let $\ell\nmid|\Gamma|$ be prime, and let
$V$ be a finite irreducible $\F_\ell[H\rt\Gamma]$-module. Set
$V^\vee:=\Hom_{\F_\ell}(V,\F_\ell)$. We call $V$ eligible at
$(\mathcal L,\boldsymbol H)$ if $V\rt H\in\mathcal L$, as in
\cite[Section~1.6]{SawinWood}. The section $H\to V\rt H$ pushes
$s_H$ forward to an orientation for which the projection is an
oriented morphism. Its kernel has no nontrivial proper subgroup
normal in the source and stable under $\Gamma$, since $V$ is irreducible.

Following \cite[Definition~7.28]{SawinWood}, write
\[
 H^2(H\rt\Gamma,V)^{\mathcal L}
\]
for the classes whose associated $\Gamma$-extensions of $H$ by $V$
have underlying group in $\mathcal L$. For eligible $V$, this set
is nonempty and is an $\F_\ell$-subspace by
\cite[Lemma~5.9]{SawinWood}. Define the module weight
\begin{equation}\label{eq:def-module-weight}
 w_{H,U}(V)
 :=\frac{\prod_{E\in U}|V^{\gamma_E}|}{|V^{H\rt\Gamma}|}
 =\frac{|V^\Gamma|}{|V^{H\rt\Gamma}|}\,w_U(V).
\end{equation}
Here $w_U(V)$ applies the earlier group weight to the additive
$\Gamma$-group of $V$. The additional subscript $H$ records the use of
the full $H\rt\Gamma$-fixed space in the denominator. Both weights are
numbers, and all fixed spaces in this formula are counted by cardinality.

When $\ell\mid n$, we want the orientation $s_H$ to lift to the
extension $1\to V\to E\xrightarrow{\pi}H\to1$: that is, we want
$s_E\in H_3(E,\Z/n)^\Gamma$ with $\pi_*(s_E)=s_H$. Write
$H^2(H\rt\Gamma,V)^{\mathcal L,s_H}$ for the subspace of classes
$\alpha\in H^2(H\rt\Gamma,V)^{\mathcal L}$ admitting such a lift.
By \cite[Definition~7.28 and Lemma~7.29]{SawinWood}, this is equivalent to
\[
\begin{aligned}
 s_H(\alpha\smile\beta)&=0
 &&\text{for every }\beta\in H^1(H\rt\Gamma,V^\vee),\\
 s_H\bigl(B(\varphi_*\alpha)\bigr)&=0
 &&\text{for every }\varphi\in\Hom_{H\rt\Gamma}(V,\Z/n).
\end{aligned}
\]
Here $s_H$ evaluates classes in $H^3(H\rt\Gamma,\Z/n)$ by
\cite[Lemmas~5.4--5.5]{SawinWood}, the cup product uses
$\Hom(V,\Z/n)\simeq V^\vee$ via $(\Z/n)[\ell]\simeq\F_\ell$,
and $B:H^2(H\rt\Gamma,\Z/n)\to H^3(H\rt\Gamma,\Z/n)$ is the Bockstein for
$0\to\Z/n\xrightarrow{\times n}\Z/n^2\to\Z/n\to0$.
The second condition is vacuous for nontrivial irreducible $V$;
for the trivial module $V=\F_\ell$, it adds at most one $\F_\ell$-linear condition.

If $H_\Gamma=1$ and $C(\mathcal L,\boldsymbol H)$ has positive
measure, \cite[Proposition~7.49]{SawinWood} gives, for each eligible $V$,
\begin{align}
 |H^2(H\rt\Gamma,V)^{\mathcal L}|
 &\le |H^1(H\rt\Gamma,V)|w_{H,U}(V),
 &&\ell\nmid n,
 \label{eq:def-support-prime}\\
 |H^1(H\rt\Gamma,V^\vee)|\,
 |H^2(H\rt\Gamma,V)^{\mathcal L,s_H}|
 &\le |H^1(H\rt\Gamma,V)|w_{H,U}(V),
 &&\ell\mid n.
 \label{eq:def-support-dividing}
\end{align}

\subsection{A uniform finite-level cohomology bound}

\begin{lemma}\label{lem:def-finite}
Let $\boldsymbol X\in\Omega_+$, let
$\boldsymbol H=\boldsymbol X^{\mathcal L}$, and let $V$ be eligible
at $(\mathcal L,\boldsymbol H)$. Then
\begin{equation}\label{eq:def-finite}
 \dim_{\F_\ell}H^2(H\rt\Gamma,V)^{\mathcal L}
 \le \dim_{\F_\ell}H^1(H\rt\Gamma,V)+|U|\dim_{\F_\ell}V.
\end{equation}
If $\ell\nmid n$, the more precise bound is
\begin{equation}\label{eq:def-finite-refined}
 \dim_{\F_\ell}H^2(H\rt\Gamma,V)^{\mathcal L}
 \le \dim_{\F_\ell}H^1(H\rt\Gamma,V)
 +\sum_{E\in U}\dim_{\F_\ell}V^{\gamma_E}
                  -\dim_{\F_\ell}V^{H\rt\Gamma}.
\end{equation}
\end{lemma}

\begin{proof}
$C(\mathcal L,\boldsymbol H)$ has positive measure by \cref{lem:def-positive}.
Throughout the proof, all dimensions are over $\F_\ell$. Since $V$ is irreducible,
\[
 \dim V^{H\rt\Gamma}
 =\begin{cases}
    0,&V\text{ is nontrivial},\\
    1,&V=\F_\ell\text{ is trivial}.
   \end{cases}
\]
Equation \eqref{eq:def-module-weight} gives
\begin{equation}\label{eq:def-fixed-spaces}
 \log_\ell w_{H,U}(V)
   =\sum_{E\in U}\dim V^{\gamma_E}-\dim V^{H\rt\Gamma}
   \le |U|\dim V-\dim V^{H\rt\Gamma}.
\end{equation}
For $\ell\nmid n$, take logarithms in
\eqref{eq:def-support-prime}. Keeping the equality in
\eqref{eq:def-fixed-spaces} gives \eqref{eq:def-finite-refined}; its
upper bound gives \eqref{eq:def-finite}.

Now suppose $\ell\mid n$.
The cup-product conditions defining $ H^2(H\rt\Gamma,V)^{\mathcal L,s_H}$ give a linear map
from $ H^2(H\rt\Gamma,V)^{\mathcal L}$ to the dual of $H^1(H\rt\Gamma,V^\vee)$. The Bockstein conditions
contribute at most $\dim V^{H\rt\Gamma}$ further independent linear
functionals. Consequently,
\[
 \dim  H^2(H\rt\Gamma,V)^{\mathcal L}-\dim  H^2(H\rt\Gamma,V)^{\mathcal L,s_H}
   \le\dim H^1(H\rt\Gamma,V^\vee)+\dim V^{H\rt\Gamma}.
\]
On the other hand, \eqref{eq:def-support-dividing} gives
\[
 \dim  H^2(H\rt\Gamma,V)^{\mathcal L,s_H}\le\dim H^1(H\rt\Gamma,V)
       +\log_\ell w_{H,U}(V)
       -\dim H^1(H\rt\Gamma,V^\vee).
\]
Adding these inequalities and using \eqref{eq:def-fixed-spaces}
proves \eqref{eq:def-finite}. In the trivial-module case, the extra
Bockstein condition is canceled by the denominator $|V^{H\rt\Gamma}|$
in \eqref{eq:def-module-weight}.
\end{proof}

\subsection{Passage to the profinite group}

The preceding bound concerns only extension classes lying in a
specified level. To apply the relation formula, we need the same bound
for all continuous second-cohomology classes.

\begin{proposition}\label{prop:def-limit}
Let $\boldsymbol X=(X,s_X)\in\Omega_+$, and let $V$ be a finite
irreducible continuous $\F_\ell[X\rt\Gamma]$-module, where
$\ell\nmid|\Gamma|$. If $H^1(X\rt\Gamma,V)$ is finite-dimensional,
then
\begin{equation}\label{eq:def-limit}
 \dim_{\F_\ell}H^2(X\rt\Gamma,V)
 \le\dim_{\F_\ell}H^1(X\rt\Gamma,V)+|U|\dim_{\F_\ell}V.
\end{equation}
If $\ell\nmid n$, then moreover
\begin{equation}\label{eq:def-limit-refined}
 \dim_{\F_\ell}H^2(X\rt\Gamma,V)
 \le\dim_{\F_\ell}H^1(X\rt\Gamma,V)
 +\sum_{E\in U}\dim_{\F_\ell}V^{\gamma_E}
                  -\dim_{\F_\ell}V^{X\rt\Gamma}.
\end{equation}
\end{proposition}

\begin{proof}
Let $ T\subseteq H^2(X\rt\Gamma,V)$ be any finite-dimensional
subspace, with basis $\alpha_1,\ldots,\alpha_a$. Since continuous
cochains with finite coefficients factor through finite quotients,
there is an open normal $\Gamma$-stable subgroup $N$ of $X$, acting
trivially on $V$, such that these classes inflate from
$H^2((X/N)\rt\Gamma,V)$. By \cite[Lemma~5.8]{SawinWood}, the lifted
classes correspond to finite $\Gamma$-extensions
$1\to V\to E_i\to X/N\to1$.

Let $\mathcal L$ be the level generated by $V\rt(X/N)$ and the
$E_i$, and put $\boldsymbol H=\boldsymbol X^{\mathcal L}$.
The quotient $X/N$ belongs to $\mathcal L$, so the map
$X\twoheadrightarrow X/N$ factors through $H$. The pullbacks
\[
 H\times_{X/N}E_i,
 \qquad
 H\times_{X/N}(V\rt(X/N))\simeq V\rt H
\]
are subdirect products of groups in $\mathcal L$, and hence belong
to $\mathcal L$. Thus $V$ is eligible at $(\mathcal L,\boldsymbol H)$,
and the classes inflated to $H\rt\Gamma$ lie in
$H^2(H\rt\Gamma,V)^{\mathcal L}$. They are linearly independent
because their further inflations to $X\rt\Gamma$ form the chosen basis of $ T$.
Applying \cref{lem:def-finite} and the inflation--restriction injection
$H^1(H\rt\Gamma,V)\hookrightarrow H^1(X\rt\Gamma,V)$ gives
\[
 \dim  T
 \le \dim H^2(H\rt\Gamma,V)^{\mathcal L}
 \le \dim H^1(X\rt\Gamma,V)+|U|\dim V.
\]
This bounds every finite-dimensional subspace $ T$ and proves
\eqref{eq:def-limit}. If $\ell\nmid n$, use
\eqref{eq:def-finite-refined} in place of the coarse bound. The fixed
spaces are unchanged, and $V^{H\rt\Gamma}=V^{X\rt\Gamma}$,
because the action factors through $H\rt\Gamma$. The same argument
then proves \eqref{eq:def-limit-refined}.
\end{proof}

\subsection{Relation ranks and consequences}

For a finitely $\Gamma$-generated pro-$|\Gamma|'$ group $G$ and a
finite irreducible $\F_\ell[G\rt\Gamma]$-module $V$, put
\[
 \xi_G(V):=\dim_{\F_\ell}\bigl(V^\Gamma/V^{G\rt\Gamma}\bigr)
 \ge0.
\]
Let $F=F_{d_\Gamma(G)}'(\Gamma)$, and let $R$ be the kernel of a
$\Gamma$-equivariant epimorphism $F\twoheadrightarrow G$. We use the
relative relation formula under the condition
$R^{\mathrm{ab}}:=R/\overline{[R,R]}\ne1$.
The relation-module calculation in
\cite[Propositions~3.4 and~3.7]{LiuPresentations} gives the generator
rank of $R^{\mathrm{ab}}$ as a profinite $G\rt\Gamma$-module.
The normal-generation argument of
\cite[Theorem~3.1]{LubotzkyPresentations}
(cf.\ \cite[Theorem~5.1]{LubotzkyPresentations}) applies with ambient
action $F\rt\Gamma$: the nonabelian simple factors together require
at most one normal generator, which is absorbed by the nonzero abelian
relation module. Thus, under this condition,
\begin{equation}\label{eq:def-Liu}
r_\Gamma(G)=d_\Gamma(G)+
\sup_{\substack{\ell\nmid|\Gamma|\\ V\ \mathrm{irreducible}}}
\left\lceil
\frac{\dim_{\F_\ell}H^2(G\rt\Gamma,V)-\dim_{\F_\ell}H^1(G\rt\Gamma,V)-\xi_G(V)}
{\dim_{\F_\ell}V}
\right\rceil.
\end{equation}
Here $V$ ranges over finite continuous modules.
We will apply the formula only after proving that all the cohomology
groups in it are finite-dimensional.

\begin{proof}[Proof of \cref{thm:bounded-deficiency}]
Intersect $\Omega_+$ with the measurable set
$\{\boldsymbol X:d_\Gamma(X)<\infty\}$, which has measure one by
\cref{thm:sharp-tail}. Fix a point in this intersection. If $d_\Gamma(X)=0$,
then $X=1$ and the conclusion is immediate.
Assume $d_\Gamma(X)>0$.

Every finite $\Gamma$-quotient of $X$ factors through a finite-level
quotient, so \cref{lem:def-positive} implies $X_\Gamma=1$.
Put $F=F_{d_\Gamma(X)}'(\Gamma)$ and $R=\ker(F\twoheadrightarrow X)$.
For a prime $\ell\nmid|\Gamma|$, freeness and $d_\Gamma(X)>0$ give
a nonzero $\Gamma$-equivariant homomorphism $f:F\to\F_\ell$, with
trivial $\Gamma$-action on the target. If $f(R)=0$, then $f$ would
factor through $X_\Gamma=1$, a contradiction. Hence
$R^{\mathrm{ab}}\ne1$, verifying the condition for \eqref{eq:def-Liu}.

The same argument applies to any finitely $\Gamma$-generated
$H\ne1$ with $H_\Gamma=1$: for every $m\ge d_\Gamma(H)$,
the kernel of $F_m'(\Gamma)\twoheadrightarrow H$ has nontrivial
abelianization. Applying the relation-module formula
\cite[Proposition~3.7]{LiuPresentations} and the preceding
normal-generation argument gives
\[
 r_{\Gamma,m}(H)-m=r_\Gamma(H)-d_\Gamma(H).
\]
For $H=1$ this holds as well, since $r_{\Gamma,m}(1)=m$.
This justifies the independence claimed in the introduction.

Let $V$ be a finite irreducible $\F_\ell[X\rt\Gamma]$-module with
$\ell\nmid|\Gamma|$. Since $H^1(\Gamma,V)=0$, every $1$-cocycle
is cohomologous to one vanishing on $\Gamma$. Such a cocycle is
determined by its values on any $d_\Gamma(X)$ $\Gamma$-generators of $X$.
It follows that
\[
 \dim_{\F_\ell}H^1(X\rt\Gamma,V)
   \le d_\Gamma(X)\dim_{\F_\ell}V<\infty.
\]
By \cref{prop:def-limit},
\[
 \dim H^2(X\rt\Gamma,V)-\dim H^1(X\rt\Gamma,V)
    \le |U|\dim V.
\]
Since $\xi_X(V)\ge0$ and $|U|$ is an integer, every ceiling in
\eqref{eq:def-Liu} is at most $|U|$. This proves
$r_\Gamma(X)\le d_\Gamma(X)+|U|$ and
$\operatorname{def}_\Gamma(X)\ge-|U|$ on a set of measure one.
\end{proof}

\section{Almost-sure FAb under a representation-theoretic hypothesis}\label{sec:fab}

We prove the sufficient condition in \cref{thm:fab}. Recall that a
profinite group is \emph{FAb} if every open subgroup has finite
abelianization. All abelianizations and normal closures in this
section are topological. We first reduce the theorem to two
probability estimates, state those estimates, and then prove them
in separate subsections.

\begin{lemma}
\label{lem:fab-normal-reduction}
Let $X$ be a profinite $\Gamma$-group. If $H^{\ab}$ is finite for
every open normal $\Gamma$-stable subgroup $H$ of $X$, then $X$ is FAb.
\end{lemma}

\begin{proof}
Let $J\le X$ be open. Its core in $X\rt\Gamma$,
\[
 H=\bigcap_{g\in X\rt\Gamma}gJg^{-1},
\]
is an open normal $\Gamma$-stable subgroup of $X$ contained in $J$.
Indeed, $J$ has finite index in $X\rt\Gamma$, so only finitely many
conjugates occur in the intersection. The image of $H$ in
$J^{\ab}$ is a quotient of the finite group $H^{\ab}$. The quotient
of $J^{\ab}$ by that image is a quotient of the finite group $J/H$.
Thus $J^{\ab}$ is finite, as required.
\end{proof}

Throughout this section, $H$ denotes an open normal $\Gamma$-stable
subgroup of $X$, and, when a name for the quotient map is needed, we write
\[
 f:X\onto X/H.
\]
Conjugation and the $\Gamma$-action give $H^{\ab}$ an action of
$(X/H)\rt\Gamma$: inner automorphisms of $H$ act trivially on its
abelianization, so the action of $X\rt\Gamma$ factors through
\[
 (X\rt\Gamma)/H\simeq (X/H)\rt\Gamma.
\]
Write $(H^{\ab})_p$ for the pro-$p$ primary component of $H^{\ab}$.

\begin{samepage}
\begin{lemma}
\label{lem:fab-primary-reduction}
Suppose the following two assertions hold:
\begin{enumerate}
\item For every prime $p$, with probability one, $(H^{\ab})_p$ is
finite for every open normal $\Gamma$-stable subgroup $H$ of $X$.
\item With probability one, for every such $H$ there is a finite set
of primes $S_H$ such that $(H^{\ab})_p=0$ for every $p\notin S_H$.
\end{enumerate}
Then $\Prob(X\text{ is FAb})=1$.
\end{lemma}
\end{samepage}

\begin{proof}
Intersect the probability-one events in the first assertion over
all primes, and then intersect with the event in the second assertion.
On the resulting event, for every open normal $\Gamma$-stable $H$,
\[
 H^{\ab}\simeq\prod_{p\in S_H}(H^{\ab})_p
\]
is a finite product of finite groups. Apply
\cref{lem:fab-normal-reduction}.
\end{proof}

We want to express the two finiteness conditions in
\cref{lem:fab-primary-reduction} in terms of finite extensions and
cohomology. First, we need a way to range over all open normal
$\Gamma$-stable subgroups of $X$. These are precisely the kernels
of continuous $\Gamma$-equivariant epimorphisms onto finite
$\Gamma$-groups. We therefore fix a finite $\Gamma$-group $P$ and
consider all maps
\[
 f\in\GammaSur(X,P),\qquad H=\ker f,\qquad X/H\simeq_\Gamma P.
\]
We will prove the two conditions simultaneously for all these
subgroups $H$. Since there are only countably many finite
$\Gamma$-groups up to isomorphism, we can then intersect the
resulting probability-one events over $P$.

\medskip\noindent\textbf{The first condition: excluding free $\Z_p$-quotients.}
Fix a prime $p$. Almost surely, $X$ is topologically finitely
generated as an ordinary profinite group by \cref{thm:sharp-tail}.
Since $H$ has finite index in
$X$, the profinite Schreier theorem makes $H$, and hence
$H^{\ab}$, topologically finitely generated. Thus $(H^{\ab})_p$
is a finitely generated $\Z_p$-module. It is finite unless its
quotient by its finite torsion subgroup is a nonzero finite-rank
free $\Z_p$-module. The torsion subgroup is preserved by the
$P\rt\Gamma$-action, so this quotient would be a
nonzero $\Z_p[P\rt\Gamma]$-module
$\Lambda$ that is free of finite rank over $\Z_p$. It is therefore
enough to exclude every nonzero equivariant quotient of $H^{\ab}$
that is free of finite rank as a $\Z_p$-module.
The transfer argument in \cref{lem:fab-transfer} shows that any
such quotient satisfies $\Lambda^{P\rt\Gamma}=0$.

For $p\nmid|\Gamma|$, fix such a module $\Lambda$ and put
\[
 \Lambda_t=\Lambda/p^t\Lambda,\qquad
 \upsilon_\Lambda=
   \sum_{E\in U}\operatorname{rank}_{\Z_p}\Lambda^{\gamma_E}.
\]
An equivariant quotient $H^{\ab}\onto\Lambda$ would give, for every
$t\ge1$, a finite $\Gamma$-quotient of $X$ in an extension
\begin{equation}\label{eq:fab-extension-family}
 1\longrightarrow\Lambda_t\longrightarrow\widetilde X_t
       \longrightarrow P\longrightarrow1.
\end{equation}
Indeed, the kernel of $H\onto\Lambda_t$ is normal in $X$ and
$\Gamma$-stable. To rule out such a quotient, it is enough to
show that the probability of having \emph{any} extension
\eqref{eq:fab-extension-family} as a quotient tends to zero with $t$.

These extensions are classified by the group
$H^2(P\rt\Gamma,\Lambda_t)$; see \cite[Lemma~5.8]{SawinWood}.
For each $\alpha\in H^2(P\rt\Gamma,\Lambda_t)$, let
$\widetilde X_{t,\alpha}$ be a representative extension.
The union bound and \cref{prop:underlying-moment} give
\begin{equation}\label{eq:fab-fixed-cohomological-reduction}
\begin{aligned}
 \Prob\bigl(\exists\,\alpha\in H^2(P\rt\Gamma,\Lambda_t):
                   X\onto_\Gamma\widetilde X_{t,\alpha}\bigr)
    &\le\sum_{\alpha\in H^2(P\rt\Gamma,\Lambda_t)}
       \Prob(X\onto_\Gamma\widetilde X_{t,\alpha})\\
       &\le\sum_{\alpha\in H^2(P\rt\Gamma,\Lambda_t)}
 \frac{|H^2(\widetilde X_{t,\alpha}\rt\Gamma,\Z/n)|}
      {w_U(\widetilde X_{t,\alpha})\,
       |\GammaAut(\widetilde X_{t,\alpha})|}.
\end{aligned}
\end{equation}
The sum may count a quotient type more than once, which only
increases the bound. A target with nontrivial $\Gamma$-coinvariants
has zero quotient probability, so including its nonnegative term
also causes no problem. We have now expressed the required
probability bound in cohomological terms. The estimate to prove is
the following.

\begin{lemma}
\label{lem:fab-primary-kernels}
Suppose that every nontrivial irreducible complex
$\Gamma$-module $W$ satisfies
\[
 \sum_{E\in U}\dim_{\mathbb C}W^{\gamma_E}>0.
\]
Fix a finite $\Gamma$-group $P$ as above, a prime
$p\nmid|\Gamma|$, and a nonzero $\Z_p[P\rt\Gamma]$-module $\Lambda$ that
is free of finite rank as a $\Z_p$-module and satisfies
$\Lambda^{P\rt\Gamma}=0$. Then $\upsilon_\Lambda>0$ and
\begin{equation}\label{eq:fab-fixed-prime-target}
 \sum_{\alpha\in H^2(P\rt\Gamma,\Lambda_t)}
 \frac{|H^2(\widetilde X_{t,\alpha}\rt\Gamma,\Z/n)|}
      {w_U(\widetilde X_{t,\alpha})\,
       |\GammaAut(\widetilde X_{t,\alpha})|}
       \ll_{\Gamma,n,U,P,p,\Lambda}p^{-t\upsilon_\Lambda}
       \qquad(t\ge1).
\end{equation}
The implied constant is independent of $t$ and of the extension
class. Consequently, for every prime $p$, with probability one,
$(H^{\ab})_p$ is finite for every open normal $\Gamma$-stable
subgroup $H$ of $X$.
\end{lemma}

We prove \eqref{eq:fab-fixed-prime-target} in
\cref{subsec:fab-fixed-prime} by bounding the number of extension classes in
$H^2(P\rt\Gamma,\Lambda_t)$, the cohomology numerator of each summand,
and the reciprocal weight--automorphism factor. The first two
bounds are independent of $t$; the third supplies the decay
$p^{-t\upsilon_\Lambda}$. This excludes each fixed equivariant quotient with a free
finite-rank $\Z_p$-module as target. Countability of the possible
isomorphism classes of these targets then excludes all of them.

\medskip\noindent\textbf{The second condition: a summable bound across primes.}
For the same fixed $P$, we need to bound, by a quantity summable
over $p$, the probability that there exists $f\in\GammaSur(X,P)$
with $((\ker f)^{\ab})_p\ne0$. If such a map exists, writing
$H=\ker f$, the nonzero finitely generated module $(H^{\ab})_p$
has a nonzero finite reduction modulo $p$.
An irreducible $\F_p[P\rt\Gamma]$-quotient $V$ of that reduction
produces a $\Gamma$-equivariant extension of $P$ by $V$ which is a
quotient of $X$. Thus, as in the first condition,
we need to bound the probability of having any one of a family of
finite extensions as a quotient.

It is enough to treat $p\nmid n|P|\,|\Gamma|$, since only finitely
many primes are excluded. For these primes,
$H^2(P\rt\Gamma,V)=0$, so every such extension is
$\Gamma$-isomorphic to $V\rt P$.
Moreover, $V$ cannot be the trivial $P\rt\Gamma$-module: otherwise
$V\times P$ would have a nontrivial quotient with trivial
$\Gamma$-action, contrary to $X_\Gamma=1$.

Write $\operatorname{Irr}_p^*(P\rt\Gamma)$ for the isomorphism
classes of nontrivial irreducible $\F_p[P\rt\Gamma]$-modules.
Applying the same quotient-moment
bound now gives
\begin{equation}\label{eq:fab-prime-cohomological-reduction}
\begin{aligned}
 \Prob\bigl(\exists f\in\GammaSur(X,P):
          ((\ker f)^{\ab})_p\ne0\bigr)&\le \sum_{V\in\operatorname{Irr}_p^*(P\rt\Gamma)}
          \Prob(X\onto_\Gamma (V\rt P))\\
 &\le \sum_{V\in\operatorname{Irr}_p^*(P\rt\Gamma)}
 \frac{|H^2((V\rt P)\rt\Gamma,\Z/n)|}
      {w_U(V\rt P)\,|\GammaAut(V\rt P)|}.
\end{aligned}
\end{equation}
Thus the second condition will follow from the following
cohomological bound.

\begin{lemma}
\label{lem:fab-large-primes}
Suppose that every nontrivial irreducible complex
$\Gamma$-module $W$ satisfies
\[
 \sum_{E\in U}\dim_{\mathbb C}W^{\gamma_E}>0.
\]
For every fixed finite $\Gamma$-group $P$ with
$P_\Gamma=1$ and every prime $p\nmid n|P|\,|\Gamma|$,
\begin{equation}\label{eq:fab-large-prime-estimate}
 \sum_{V\in\operatorname{Irr}_p^*(P\rt\Gamma)}
 \frac{|H^2((V\rt P)\rt\Gamma,\Z/n)|}
      {w_U(V\rt P)\,|\GammaAut(V\rt P)|}
 \le \frac{|P|\,|\Gamma|}{p(p-1)}
       \frac{|H^2(P\rt\Gamma,\Z/n)|}{w_U(P)}.
\end{equation}
The factor $|P|\,|\Gamma|\,|H^2(P\rt\Gamma,\Z/n)|/w_U(P)$
is finite and independent of $p$. Consequently, for this fixed
$P$, almost surely a single finite set of primes contains every
$p$ for which $(H^{\ab})_p\ne0$, simultaneously for all $H=\ker f$
with $f\in\GammaSur(X,P)$. In particular, almost surely every open
normal $\Gamma$-stable subgroup $H$ of $X$ has only finitely many
nonzero primary components in its abelianization.
\end{lemma}

We prove this bound in \cref{subsec:fab-across-primes}. Its three
factors are the number of possible irreducible modules, the fixed
cohomology factor $|H^2(P\rt\Gamma,\Z/n)|/w_U(P)$, and a
weight--automorphism ratio bounded by $1/[p(p-1)]$. Combining the
bound with \eqref{eq:fab-prime-cohomological-reduction} and
Borel--Cantelli gives the assertion for each $P$. A countable
intersection over $P$ gives the simultaneous conclusion required
by \cref{lem:fab-primary-reduction}.

\subsection{Bounds at a fixed prime}
\label{subsec:fab-fixed-prime}

We prove \cref{lem:fab-primary-kernels} by isolating three factors
in the quotient-moment bound and estimating them separately.
We first record the representation-theoretic and transfer facts
needed for those estimates and for their application to $H^{\ab}$.

\begin{lemma}\label{lem:fab-character-comparison}
The assertion that every nontrivial irreducible complex
$\Gamma$-module $W$ satisfies
\[
 \sum_{E\in U}\dim_{\mathbb C}W^{\gamma_E}>0
\]
is equivalent to $u_V>0$ for every $V\in\cI$.
If these equivalent assertions hold, then, for any nonzero
finite-dimensional $\Gamma$-representation $D$ in characteristic
zero or in characteristic $p\nmid|\Gamma|$,
\[
 \sum_{E\in U}\dim D^{\gamma_E}>0.
\]
In particular, if $p\nmid|\Gamma|$ and $\Lambda$ is a nonzero
$\Z_p[\Gamma]$-module that is free of finite rank over $\Z_p$, then
\begin{equation}\label{eq:fab-free-module-weight}
 \upsilon_\Lambda>0,
 \qquad
 \prod_{E\in U}|(\Lambda/p^t\Lambda)^{\gamma_E}|
       =p^{t \upsilon_\Lambda}.
\end{equation}
If $n>1$, the Sawin--Wood fixed-vector hypothesis ensures
$\sum_{E\in U}\dim_{\mathbb C}W^{\gamma_E}>0$ for every
nontrivial irreducible complex $\Gamma$-module $W$.
\end{lemma}

\begin{proof}
Over an algebraically closed field of characteristic $p\nmid|\Gamma|$,
the ordinary and modular irreducible characters correspond, and the
multiplicity of the eigenvalue $1$ for any element is preserved; see
\cite{IsaacsBook,NavarroBook}. Equivalently, restrict to its cyclic
subgroup and compare the multiplicity of the trivial character.
Extension of scalars does not change fixed-space dimensions. Decomposing
an irreducible module over $\F_p$ over an algebraic closure
therefore proves the equivalence with $u_V>0$. The same comparison in
characteristic zero applies to representations over $\mathbb Q_p$.

All the relevant representations are semisimple. Their nontrivial
constituents satisfy the required positivity by the assumed
inequality; their trivial constituents satisfy it because
$U\ne\varnothing$. This proves
the assertion for $D$. Apply it to $\Lambda\otimes\mathbb Q_p$ to get
$\upsilon_\Lambda>0$. The order of each $\gamma_E$ is a unit in $\Z_p$,
so averaging makes invariants exact and gives
\[
 (\Lambda/p^t\Lambda)^{\gamma_E}
       =\Lambda^{\gamma_E}/p^t\Lambda^{\gamma_E}.
\]
This proves \eqref{eq:fab-free-module-weight}. Finally, if $p\mid n$,
the Sawin--Wood fixed-vector hypothesis gives $u_V>0$ for
every irreducible $\F_p[\Gamma]$-module; character comparison
then gives the asserted inequality for every nontrivial
irreducible complex $\Gamma$-module.
\end{proof}

\begin{lemma}\label{lem:fab-free-module-types}
For a finite group $J$, a prime $p$, and an integer $d\ge0$,
there are only finitely many isomorphism classes of $\Z_p[J]$-modules
that are free of rank $d$ as $\Z_p$-modules.
\end{lemma}

\begin{proof}
Put $c=v_p(|J|)$. Suppose two such modules have isomorphic
reductions modulo $p^{c+1}$, and lift such an isomorphism to an invertible
$\Z_p$-linear map $T$. Write $\rho_1,\rho_2$ for their actions.
Every difference $\rho_2(g)T\rho_1(g)^{-1}-T$ is divisible by
$p^{c+1}$. Hence
\[
 T_0=\frac1{|J|}\sum_{g\in J}\rho_2(g)T\rho_1(g)^{-1}
\]
is integral, is equivariant, and is congruent to $T$ modulo $p$.
It is therefore an isomorphism. There are only finitely many actions
on $(\Z/p^{c+1})^d$ for fixed $d$, which proves the assertion.
No assumption $p\nmid|J|$ is needed here.
\end{proof}

\begin{lemma}\label{lem:fab-transfer}
Suppose $X_\Gamma=1$ and $p\nmid|\Gamma|$. If a nonzero
$\Z_p[(X/H)\rt\Gamma]$-module $\Lambda$ that is free of finite
rank over $\Z_p$ is an $(X/H)\rt\Gamma$-equivariant quotient
of $H^{\ab}$, then
\[
 (\Lambda\otimes\mathbb Q_p)^{(X/H)\rt\Gamma}=0.
\]
\end{lemma}

\begin{proof}
Put $\Pi=X\rt\Gamma$ and let
$\varphi:H^{\ab}\onto\Lambda$ be the given equivariant map.
We use the continuous transfer homomorphism
\[
 \operatorname{Ver}_{\Pi/H}:\Pi\longrightarrow H^{\ab}.
\]
Since $H$ is normal, the restriction of transfer to $H$ is the
norm operator:
\[
 \operatorname{Ver}_{\Pi/H}(h)
   =\sum_{g\in(X/H)\rt\Gamma}g\cdot[h]
       \qquad(h\in H).
\]

The condition $X_\Gamma=1$ gives
$\Pi^{\ab}\simeq\Gamma^{\ab}$. Thus the composition
\[
 \Pi\xrightarrow{\operatorname{Ver}_{\Pi/H}}
 H^{\ab}\xrightarrow{\varphi}\Lambda
\]
is zero: it factors through the finite group $\Gamma^{\ab}$, whereas
$\Lambda$ is torsion-free. Restricting to $H$ and using the
surjectivity of $\varphi$ gives
\[
 \left(\sum_{g\in(X/H)\rt\Gamma}g\right)\Lambda=0.
\]
Over $\mathbb Q_p$, the operator
\[
 \frac1{|X/H|\,|\Gamma|}\sum_{g\in(X/H)\rt\Gamma}g
\]
is the projection onto $(X/H)\rt\Gamma$-invariants. Hence
$(\Lambda\otimes\mathbb Q_p)^{(X/H)\rt\Gamma}=0$, as claimed.
The division takes place in $\mathbb Q_p$, so it is valid even
when $p\mid|X/H|$.
\end{proof}

\begin{proof}[Proof of \cref{lem:fab-primary-kernels}]
Fix $P$, $p$, and $\Lambda$ as in the statement. We first bound
the sum over extension classes in
\eqref{eq:fab-fixed-cohomological-reduction}, and then explain why
this excludes every nonzero equivariant quotient of $H^{\ab}$
that is free of finite rank over $\Z_p$.

\medskip\noindent\textbf{1. The target estimate and the three factors to bound.}
By \eqref{eq:fab-fixed-cohomological-reduction}, the probability
that some extension is a quotient is bounded by the sum appearing
there. We need to show that this sum tends to zero with $t$,
with constants independent of $t$ and of the extension class.
The extensions $\widetilde X_{t,\alpha}$ of $P$ by $\Lambda_t$
are indexed by $\alpha\in H^2(P\rt\Gamma,\Lambda_t)$.
Taking separate maxima in the sum gives
\begin{equation}\label{eq:fab-fixed-prime-three-factors}
\begin{aligned}
&\sum_{\alpha\in H^2(P\rt\Gamma,\Lambda_t)}
 \frac{|H^2(\widetilde X_{t,\alpha}\rt\Gamma,\Z/n)|}
      {w_U(\widetilde X_{t,\alpha})\,
       |\GammaAut(\widetilde X_{t,\alpha})|}\le \\
&|H^2(P\rt\Gamma,\Lambda_t)|
\max_{\alpha\in H^2(P\rt\Gamma,\Lambda_t)}
       |H^2(\widetilde X_{t,\alpha}\rt\Gamma,\Z/n)|
       \max_{\alpha\in H^2(P\rt\Gamma,\Lambda_t)}
 \frac1{w_U(\widetilde X_{t,\alpha})\,
             |\GammaAut(\widetilde X_{t,\alpha})|}.
\end{aligned}
\end{equation}
It therefore suffices to show that the first two factors on the
right are bounded independently of $t$, while the third is
$O(p^{-t\upsilon_\Lambda})$. We prove these bounds in order.

\medskip\noindent\textbf{2. Bounding the number of extension classes.}
The coefficient sequence
\[
 0\longrightarrow\Lambda\xrightarrow{p^t}\Lambda
   \longrightarrow\Lambda_t\longrightarrow0
\]
gives an exact sequence
\[
 H^2(P\rt\Gamma,\Lambda)\longrightarrow H^2(P\rt\Gamma,\Lambda_t)
       \longrightarrow H^3(P\rt\Gamma,\Lambda).
\]
Consequently,
\[
 |H^2(P\rt\Gamma,\Lambda_t)|
 \le |H^2(P\rt\Gamma,\Lambda)|\,
       |H^3(P\rt\Gamma,\Lambda)|=O(1).
\]
The groups on the right are finite: positive-degree cohomology of
the finite group $P\rt\Gamma$ with coefficients in the finite-rank
free $\Z_p$-module $\Lambda$ is finitely generated over $\Z_p$ and killed by $|P|\,|\Gamma|$.
The number of extension classes is therefore bounded independently
of $t$.

\medskip\noindent\textbf{3. Bounding the cohomology numerator.}
For each extension, use
\[
 1\longrightarrow\Lambda_t\longrightarrow
 \widetilde X_{t,\alpha}\rt\Gamma\longrightarrow P\rt\Gamma\longrightarrow1.
\]
The degree-two terms of the Lyndon--Hochschild--Serre spectral
sequence give
\[
 |H^2(\widetilde X_{t,\alpha}\rt\Gamma,\Z/n)|
 \le |H^2(P\rt\Gamma,\Z/n)|\,
       |H^1(P\rt\Gamma,H^1(\Lambda_t,\Z/n))|
 {}\times |H^2(\Lambda_t,\Z/n)|.
\]
The first factor is fixed. For the other two, put
$d_\Lambda=\operatorname{rank}_{\Z_p}\Lambda$ and $v=v_p(n)$.
Although the exponent of $\Lambda_t$ grows with $t$, its number of
cyclic generators is always $d_\Lambda$. The fixed coefficients
$\Z/n$ therefore give the uniform bounds
\[
 |H^1(\Lambda_t,\Z/n)|\le p^{v d_\Lambda},\quad
 |H^2(\Lambda_t,\Z/n)|\le
             p^{v(d_\Lambda+\binom{d_\Lambda}{2})}.
\]
The second bound follows from the universal-coefficient sequence
for a finite abelian group. Moreover, for any finite $P\rt\Gamma$-module $T$,
\[
 |H^1(P\rt\Gamma,T)|\le |T|^{|P|\,|\Gamma|},
\]
since its cocycles are functions from $P\rt\Gamma$ to $T$. Applying this to
$T=H^1(\Lambda_t,\Z/n)$ bounds the middle spectral-sequence factor
independently of both $t$ and $\alpha$. Thus the second factor in
\eqref{eq:fab-fixed-prime-three-factors} is uniformly bounded.
This argument also covers $p\mid n$ and $p\mid|P|$; neither is
excluded in the fixed-prime estimate.

\medskip\noindent\textbf{4. Bounding the weight--automorphism factor.}
Suppress $\alpha$ and write $\widetilde X_t=\widetilde X_{t,\alpha}$. We need to show
\[
 \frac1{w_U(\widetilde X_t)\,|\GammaAut(\widetilde X_t)|}
          \ll p^{-t\upsilon_\Lambda}.
\]
First calculate the weight. Coprime fixed-point lifting gives
\begin{equation}\label{eq:fab-weight-factorization}
 w_U(\widetilde X_t)=w_U(P)
       \frac{\prod_{E\in U}|\Lambda_t^{\gamma_E}|}
            {|\Lambda_t^\Gamma|}.
\end{equation}
For example, a $\gamma_E$-fixed coset has a fixed lift because its
obstruction lies in $H^1(\langle\gamma_E\rangle,\Lambda_t)=0$;
the same argument applies to the full subgroup $\Gamma$.
By \cref{lem:fab-character-comparison},
\[
 \prod_{E\in U}|\Lambda_t^{\gamma_E}|
        =p^{t\upsilon_\Lambda},\qquad \upsilon_\Lambda>0.
\]
Substituting this into \eqref{eq:fab-weight-factorization} isolates
exactly the ratio we still need to control:
\begin{equation}\label{eq:fab-fixed-prime-cancellation-target}
 \frac1{w_U(\widetilde X_t)\,|\GammaAut(\widetilde X_t)|}
 =\frac{p^{-t\upsilon_\Lambda}}{w_U(P)}
     \frac{|\Lambda_t^\Gamma|}{|\GammaAut(\widetilde X_t)|}.
\end{equation}
Thus the remaining task is to show that
$|\Lambda_t^\Gamma|/|\GammaAut(\widetilde X_t)|$ is bounded independently of $t$.

Conjugation by elements of $\Lambda_t^\Gamma$ defines
$\Gamma$-automorphisms of $\widetilde X_t$ fixing the quotient $P$. The kernel
of this conjugation map consists of elements fixed by both $P$
and $\Gamma$, namely $\Lambda_t^{P\rt\Gamma}$. It follows that
\[
 |\GammaAut(\widetilde X_t)|\ge
       \frac{|\Lambda_t^\Gamma|}{|\Lambda_t^{P\rt\Gamma}|},
 \qquad
 \frac{|\Lambda_t^\Gamma|}{|\GammaAut(\widetilde X_t)|}
       \le |\Lambda_t^{P\rt\Gamma}|.
\]
The hypothesis $\Lambda^{P\rt\Gamma}=0$ now controls the remaining term.
The same coefficient sequence as in step~2 gives an injection
\[
 \Lambda_t^{P\rt\Gamma}\hookrightarrow H^1(P\rt\Gamma,\Lambda).
\]
In particular,
\begin{equation}\label{eq:fab-invariants-bound}
 |\Lambda_t^{P\rt\Gamma}|\le |H^1(P\rt\Gamma,\Lambda)|=O(1).
\end{equation}
Combining this with \eqref{eq:fab-fixed-prime-cancellation-target}
yields
\[
 \frac1{w_U(\widetilde X_t)\,|\GammaAut(\widetilde X_t)|}
 \le\frac{|H^1(P\rt\Gamma,\Lambda)|}{w_U(P)}p^{-t\upsilon_\Lambda}.
\]
The bound is uniform in the extension class, so the final factor
in \eqref{eq:fab-fixed-prime-three-factors} is
$O(p^{-t\upsilon_\Lambda})$.

\medskip\noindent\textbf{5. Combining the bounds and excluding free $\Z_p$-quotients.}
By \eqref{eq:fab-fixed-prime-three-factors} and the three estimates
above, the sum in \eqref{eq:fab-fixed-cohomological-reduction} is
$O(p^{-t\upsilon_\Lambda})$, which tends to zero because
$\upsilon_\Lambda>0$. This proves \eqref{eq:fab-fixed-prime-target}.
In particular, the probability that $X$ admits a $\Gamma$-quotient
isomorphic to $\widetilde X_{t,\alpha}$ for some
$\alpha\in H^2(P\rt\Gamma,\Lambda_t)$ tends to zero as $t\to\infty$.

By the reduction preceding the lemma, it remains to exclude nonzero
$P\rt\Gamma$-equivariant quotients of $H^{\ab}$ that are free
of finite rank as $\Z_p$-modules. For any fixed $\Lambda$ as above,
such a quotient would yield a $\Gamma$-quotient of $X$ of the above
form for every $t$. The probability of this is zero, since the
probabilities for the individual values of $t$ tend to zero.

Finally, \cref{lem:fab-free-module-types} gives only countably many
isomorphism classes of $\Z_p[P\rt\Gamma]$-modules that are free
of finite rank over $\Z_p$. Taking the union of the corresponding
null sets excludes every nonzero quotient of this kind,
simultaneously for all $H=\ker f$ with $f\in\GammaSur(X,P)$.
There are only countably many finite $\Gamma$-groups up to
isomorphism, and
\cref{prop:underlying-moment} excludes those with nontrivial
coinvariants. A further countable intersection therefore gives,
for almost every sample $X$, finiteness of $(H^{\ab})_p$ for every
open normal $\Gamma$-stable subgroup $H\le X$.
\end{proof}

\medskip\noindent\textbf{An unconditional fixed-prime consequence for $X$.}
For the homology argument in \cref{sec:homology}, we also need the
following unconditional statement about individual primary
components. It does not assert finiteness of the entire abelianization.

\begin{lemma}\label{lem:fab-primary-X}
Almost surely, $(X^{\ab})_p$ is finite for every prime $p$.
\end{lemma}

\begin{proof}
For $p\mid|\Gamma|$, the primary component is zero. Otherwise,
as explained before \cref{lem:fab-primary-kernels}, if
$(X^{\ab})_p$ is infinite, it has a nonzero
$\Gamma$-equivariant quotient $\Lambda$ that is free of finite
rank over $\Z_p$. Put
\[
 \Lambda_t=\Lambda/p^t\Lambda,\qquad
 d_\Lambda=\operatorname{rank}_{\Z_p}\Lambda>0.
\]

Most of the fixed-prime argument in the proof of
\cref{lem:fab-primary-kernels} applies with $P=1$. There is only
one extension class, since $H^2(\Gamma,\Lambda_t)=0$, and the
same cohomology bound gives, with $v=v_p(n)$,
\[
 |H^2(\Lambda_t\rt\Gamma,\Z/n)|
 \le p^{v(d_\Lambda+\binom{d_\Lambda}{2})},
\]
independently of $t$. The difference is the weight--automorphism
estimate. Here we do not require that
\[
 \sum_{E\in U}\dim_{\mathbb C}W^{\gamma_E}>0
\]
for every nontrivial irreducible complex $\Gamma$-module $W$.
Instead, we use $w_U(\Lambda_t)\ge1$ (because $U\ne\varnothing$)
and the faithful action of $(\Z/p^t)^\times$ by scalar
$\Gamma$-automorphisms of $\Lambda_t$. Hence
\[
 \Prob(X\onto_\Gamma\Lambda_t)
 \le\frac{|H^2(\Lambda_t\rt\Gamma,\Z/n)|}
          {w_U(\Lambda_t)\,|\GammaAut(\Lambda_t)|}
 \le\frac{p^{v(d_\Lambda+\binom{d_\Lambda}{2})}}{\varphi(p^t)}
 \longrightarrow0.
\]
A quotient $X^{\ab}\onto_\Gamma\Lambda$ would induce such a
quotient for every $t$, and therefore occurs with probability zero.
By \cref{lem:fab-free-module-types}, only countably many isomorphism
classes of finite-rank free $\Z_p$-modules with a $\Gamma$-action
need be excluded; taking the countable intersection over primes
completes the proof.
\end{proof}

\subsection{Bounds across primes}
\label{subsec:fab-across-primes}

We now prove \cref{lem:fab-large-primes}. The fixed-prime estimate
tends to zero with $t$; the estimate here must instead be summable
as $p$ varies. We again separate the number of targets, their
cohomology factor, and the weight--automorphism factor.

\begin{proof}[Proof of \cref{lem:fab-large-primes}]
Keep $P$ fixed, and let $p\nmid n|P|\,|\Gamma|$. We bound the
sum over irreducible modules in
\eqref{eq:fab-prime-cohomological-reduction}; it bounds the
probability that some $f\in\GammaSur(X,P)$ has
$((\ker f)^{\ab})_p\ne0$. As in the preceding subsection,
we first isolate three factors and then bound them separately.
Here $p$ varies, so all constants must be independent of $p$.

\medskip\noindent\textbf{1. The target estimate and the three factors to bound.}
The reduction preceding the lemma shows why these are the relevant
quotients: if there is $f\in\GammaSur(X,P)$, with $H=\ker f$,
such that $(H^{\ab})_p\ne0$, then a nontrivial irreducible
quotient $V$ of $(H^{\ab})_p/p(H^{\ab})_p$ yields a quotient
$V\rt P$ of $X$. The vanishing $H^2(P\rt\Gamma,V)=0$ ensures
that every such extension splits, and the trivial module is
excluded by $X_\Gamma=1$.
The bound \eqref{eq:fab-prime-cohomological-reduction} thus applies.
If $\operatorname{Irr}_p^*(P\rt\Gamma)$ is empty, the sum and the
probability are both zero; henceforth suppose it is nonempty.

We next express each summand using the cohomology of the fixed
target $P$. Since $p\nmid n$, inflation gives an isomorphism
\[
 H^2(P\rt\Gamma,\Z/n)
       \simeq H^2((V\rt P)\rt\Gamma,\Z/n),
\]
because the positive-degree cohomology of the normal $p$-group $V$
with coefficients $\Z/n$ vanishes. The weight factorization
\eqref{eq:fab-weight-factorization}, with $V$ in place of
$\Lambda_t$, therefore gives the exact identity
\begin{equation}\label{eq:fab-prime-relative-cohomology}
 \frac{|H^2((V\rt P)\rt\Gamma,\Z/n)|}
       {w_U(V\rt P)\,|\GammaAut(V\rt P)|}
 =
 \frac{|H^2(P\rt\Gamma,\Z/n)|}{w_U(P)}
 \frac{|V^\Gamma|}
      {|\GammaAut(V\rt P)|\prod_{E\in U}|V^{\gamma_E}|}.
\end{equation}
This equality holds for every module in the sum, even when the
corresponding quotient event has probability zero. Summing it
gives the three-factor estimate
\begin{equation}\label{eq:fab-prime-three-factors}
\begin{aligned}
&\sum_{V\in\operatorname{Irr}_p^*(P\rt\Gamma)}
 \frac{|H^2((V\rt P)\rt\Gamma,\Z/n)|}
      {w_U(V\rt P)\,|\GammaAut(V\rt P)|}\\
&\quad\le |\operatorname{Irr}_p^*(P\rt\Gamma)|\,
 \frac{|H^2(P\rt\Gamma,\Z/n)|}{w_U(P)}\\
&\qquad{}\times\max_{V\in\operatorname{Irr}_p^*(P\rt\Gamma)}
 \frac{|V^\Gamma|}
      {|\GammaAut(V\rt P)|\prod_{E\in U}|V^{\gamma_E}|}.
\end{aligned}
\end{equation}
It suffices to show that the number of possible irreducible modules
is at most $|P|\,|\Gamma|$, that the cohomology factor is finite
and independent of $p$, and that the final maximum is at most
$1/[p(p-1)]$. The first two factors are uniformly bounded,
while the third supplies the summable decay.

\medskip\noindent\textbf{2. Bounding the number of irreducible modules.}
Because $p\nmid|P|\,|\Gamma|$, Maschke's theorem makes
$\F_p[P\rt\Gamma]$ semisimple. Each isomorphism class of simple
modules corresponds to a simple component of this algebra, and
each component has positive $\F_p$-dimension. Since the algebra
has dimension $|P|\,|\Gamma|$,
\[
 |\operatorname{Irr}_p^*(P\rt\Gamma)|
       \le |P|\,|\Gamma|.
\]
In particular, the number of possible modules is bounded independently
of the prime.

\medskip\noindent\textbf{3. Bounding the cohomology factor of the fixed quotient.}
The factor
\[
 \frac{|H^2(P\rt\Gamma,\Z/n)|}{w_U(P)}
       =\Exp|\GammaSur(X,P)|<\infty
\]
is fixed, where the equality with the moment uses $P_\Gamma=1$.
Everything here depends only on $\Gamma,n,U,P$, not on $p$ or
$V$. Thus the first two factors are uniformly bounded.

\medskip\noindent\textbf{4. Bounding the weight--automorphism factor.}
Fix $V\in\operatorname{Irr}_p^*(P\rt\Gamma)$. We must prove
\begin{equation}\label{eq:fab-prime-factor-target}
 \frac{|V^\Gamma|}
      {|\GammaAut(V\rt P)|\prod_{E\in U}|V^{\gamma_E}|}
       \le\frac1{p(p-1)}.
\end{equation}
We obtain the two denominator factors separately: automorphisms
supply $p-1$, and fixed vectors supply $p$.

\smallskip\noindent\emph{The automorphism factor supplies $p-1$.}
Write $q=|\End_{P\rt\Gamma}(V)|$ in this proof only. Schur's
lemma and the finiteness of $V$ make this endomorphism ring a finite
field containing $\F_p$, so $q\ge p$.

For each $v\in V^\Gamma$, conjugation $c_v$ is a
$\Gamma$-automorphism of $V\rt P$ fixing its quotient
$P$. Its kernel, as a map from $V^\Gamma$, is
$V^{P\rt\Gamma}=0$, since $V$ is a nontrivial irreducible
$P\rt\Gamma$-module. These inner automorphisms therefore form
a subgroup of order $|V^\Gamma|$.

Every $\lambda\in\End_{P\rt\Gamma}(V)^\times$ also defines a
$\Gamma$-automorphism
\[
 \operatorname{scal}_\lambda : V\rt P\longrightarrow V\rt P,\qquad
 \operatorname{scal}_\lambda(v,g)=(\lambda v,g)
       \quad(v\in V,\ g\in P).
\]
These form a subgroup of order $q-1$. Its intersection with the
inner subgroup is trivial, because conjugation by $V$ acts
identically on $V$, whereas $\operatorname{scal}_\lambda$ acts
there as $\lambda$. Their product therefore contains
$|V^\Gamma|(q-1)$ distinct $\Gamma$-automorphisms. Consequently,
\begin{equation}\label{eq:fab-prime-automorphisms}
 |\GammaAut(V\rt P)|\ge |V^\Gamma|(q-1)
             \ge |V^\Gamma|(p-1).
\end{equation}
In particular, the factor $|V^\Gamma|$ in
\eqref{eq:fab-prime-factor-target} cancels, leaving at most
$1/(p-1)$ from the automorphism ratio.

\smallskip\noindent\emph{The fixed-point weight supplies $p$.}
Since $V|_\Gamma\ne0$, \cref{lem:fab-character-comparison} gives
\begin{equation}\label{eq:fab-prime-fixed-vectors}
 \prod_{E\in U}|V^{\gamma_E}|
 =p^{\sum_{E\in U}\dim_{\F_p}V^{\gamma_E}}\ge p.
\end{equation}

\smallskip\noindent\emph{Combining the two factors.}
Equations \eqref{eq:fab-prime-automorphisms} and
\eqref{eq:fab-prime-fixed-vectors} give
\[
 \frac{|V^\Gamma|}
      {|\GammaAut(V\rt P)|\prod_{E\in U}|V^{\gamma_E}|}
 \le\frac{|V^\Gamma|}{|V^\Gamma|(p-1)p}
 =\frac1{p(p-1)}.
\]
The bound holds for every $V\in\operatorname{Irr}_p^*(P\rt\Gamma)$,
so the maximum in \eqref{eq:fab-prime-three-factors} is at most
$1/[p(p-1)]$, as required.

\medskip\noindent\textbf{5. Combining the bounds and summing over primes.}
Returning to \eqref{eq:fab-prime-three-factors}, the bounds for
the number of irreducible modules, the fixed cohomology factor, and
the weight--automorphism ratio give \eqref{eq:fab-large-prime-estimate}.
Combining this with \eqref{eq:fab-prime-cohomological-reduction}
shows that the probability that there is an
$f\in\GammaSur(X,P)$ with $((\ker f)^{\ab})_p\ne0$ is at most
\[
 \frac{|P|\,|\Gamma|}{p(p-1)}
       \frac{|H^2(P\rt\Gamma,\Z/n)|}{w_U(P)}.
\]
Since
\[
 \sum_{p\nmid n|P|\,|\Gamma|}\frac1{p(p-1)}<\infty,
\]
the first Borel--Cantelli lemma implies that almost surely only
finitely many primes $p$ admit any $f\in\GammaSur(X,P)$ with
$((\ker f)^{\ab})_p\ne0$. No independence is needed. Adjoining
the finitely many primes dividing $n|P|\,|\Gamma|$ gives a single
finite set containing every prime $p$ for which $(H^{\ab})_p\ne0$,
simultaneously for all $H=\ker f$ with $f\in\GammaSur(X,P)$.
Finally, take the countable intersection
over finite $\Gamma$-groups $P$ up to isomorphism. By \cref{prop:underlying-moment}, almost
every sample has no quotient with nontrivial coinvariants. Thus,
for almost every $X$, every open normal $\Gamma$-stable subgroup
$H\le X$ has an abelianization $H^{\ab}$ with only finitely many
nonzero primary components.
\end{proof}

\subsection{Proof of the FAb theorem}
\label{subsec:fab-theorem}

\begin{proof}[Proof of \cref{thm:fab}]
The hypothesis of \cref{thm:fab}, together with
\cref{lem:fab-character-comparison}, gives $u_V>0$
for every $V\in\cI$.
For every prime $p$, \cref{lem:fab-primary-kernels} gives,
almost surely, finiteness of $(H^{\ab})_p$ for every open normal
$\Gamma$-stable $H\le X$. By \cref{lem:fab-large-primes}, almost
surely every such $H$ has only finitely many nonzero primary
components. These are the two hypotheses of
\cref{lem:fab-primary-reduction}, so $\Prob(X\text{ is FAb})=1$.
\end{proof}

\section{The structure of the second homology}
\label{sec:homology}

We prove a module-surjection theorem and deduce
\cref{thm:homology-vanishing} as a consequence.
All homology is continuous homology with compact coefficients and
trivial group action on those coefficients; cohomology has discrete
coefficients. We use compact/discrete duality and the coefficient
exact sequences for profinite groups; see \cite{RibesZalesskii}.
The second homology groups are compact abelian groups with their
natural $\Gamma$-actions, and
\begin{equation}\label{eq:homology-primary-product}
 H_2(X,\widehat{\Z})=\prod_{p\nmid|\Gamma|}H_2(X,\Z_p),
 \qquad H_2(X,\Z_p)=0\quad(p\mid|\Gamma|).
\end{equation}
The last assertion follows by averaging on the finite quotients of
$X$, or by duality with discrete $p$-primary cohomology.

\subsection{The augmentation module and the structural theorem}

Put
\[
 R_\Gamma=\prod_{p\nmid|\Gamma|}\Z_p
\]
and define the compact $R_\Gamma[\Gamma]$-module
\begin{equation}\label{eq:homology-global-augmentation}
 \mathcal K_U=
 \ker\left(
   \bigoplus_{E\in U}R_\Gamma[\Gamma/\langle\gamma_E\rangle]
       \xrightarrow{\text{sum of augmentations}}R_\Gamma
      \right).
\end{equation}
Here each occurrence in $U$ contributes a summand. For
$p\nmid|\Gamma|$, its $p$-component and its reduction are
\begin{equation}\label{eq:homology-permutation-module}
 \mathcal K_{U,p}=
 \ker\left(
   \bigoplus_{E\in U}\Z_p[\Gamma/\langle\gamma_E\rangle]
       \xrightarrow{\text{sum of augmentations}}\Z_p
      \right),\quad
 \overline{\mathcal K}_{U,p}=\mathcal K_{U,p}/p\mathcal K_{U,p}.
\end{equation}
The augmentation is surjective since $U\ne\varnothing$, and splits
$\Gamma$-equivariantly by averaging. Thus $\mathcal K_{U,p}$ is a
free $\Z_p$-module, with
\begin{equation}\label{eq:homology-permutation-rank}
 \operatorname{rank}_{\Z_p}\mathcal K_{U,p}=\sum_{E\in U}[\Gamma:\langle\gamma_E\rangle]-1,
 \qquad
 \mathcal K_U=\prod_{p\nmid|\Gamma|}\mathcal K_{U,p}.
\end{equation}

For the oriented sample $\boldsymbol X=(X,s_X)$, let
\begin{equation}\label{eq:homology-orientation-boundary}
 \omega_X=\partial_n(s_X)\in H_2(X,\widehat{\Z})^\Gamma[n],
\end{equation}
where $\partial_n$ is the homological boundary associated with
\[
 0\longrightarrow\widehat{\Z}\xrightarrow{\times n}\widehat{\Z}
   \longrightarrow\Z/n\longrightarrow0.
\]
Thus $\langle\omega_X\rangle$ is finite cyclic, is fixed pointwise by
$\Gamma$, and has order dividing $n$. Its $p$-component is $\langle\omega_{X,p}\rangle$, where
$\omega_{X,p}$ is the $p$-component of $\omega_X$.

\begin{theorem}
\label{thm:homology-structure}
Under the standing assumptions on $\Gamma,n,U$, there is a measurable
set of oriented samples of $\nu_{\Gamma,n,U}$-measure one on which
there exists a continuous $R_\Gamma[\Gamma]$-linear epimorphism
\begin{equation}\label{eq:homology-structure-quotient}
 \phi_X:\mathcal K_U\onto
    H_2(X,\widehat{\Z})/\langle\omega_X\rangle.
\end{equation}
It can be lifted to a continuous $R_\Gamma[\Gamma]$-linear map
$\psi_X:\mathcal K_U\to H_2(X,\widehat{\Z})$ such that
\begin{equation}\label{eq:homology-structure-plus}
 \mathcal K_U\oplus\Z/n\onto H_2(X,\widehat{\Z}),\quad
 (x,z)\longmapsto\psi_X(x)+z\omega_X
\end{equation}
is an epimorphism. The $\Gamma$-action on $\Z/n$ is trivial, and
$1$ in that summand maps to $\omega_X$. In particular,
$H_2(X,\widehat{\Z})$ is topologically finitely generated as a profinite
abelian group. No canonical or measurable choice of the epimorphisms
is asserted.
\end{theorem}

\subsection{A cohomological bound}

We first bound the second-cohomology classes that annihilate the
orientation boundary. Universal coefficients will then give the
multiplicity bound needed for \cref{thm:homology-structure}.
For an oriented profinite $\Gamma$-group $(Y,s_Y)$, put
$\omega_Y=\partial_n(s_Y)$. If $p\nmid|\Gamma|$ and $V$ is a finite
$\F_p[\Gamma]$-module with trivial $Y$-action, define
\begin{equation}\label{eq:homology-annihilator-definition}
 \begin{split}
 H^2_{\omega_Y}(Y\rt\Gamma,V)
   :=\{\alpha\in H^2(Y\rt\Gamma,V):\;&
       \langle\alpha|_Y,\omega_Y\bmod p\rangle=0\}.
 \end{split}
\end{equation}
Here $\omega_Y\bmod p$ is its image in $H_2(Y,\F_p)$, and the
pairing takes values in $V^\Gamma$. Thus this subspace is all of
$H^2(Y\rt\Gamma,V)$ when $p\nmid n$ or $V^\Gamma=0$.
The following identity relates it to the Bockstein conditions in
\cref{subsec:positive-quotients}.

\begin{lemma}
\label{lem:homology-bockstein-pairing}
Let $J$ be a profinite group, let $s\in H_3(J,\Z/n)$, and put
$\omega=\partial_n(s)$. For $\eta\in H^2(J,\Z/n)$,
\begin{equation}\label{eq:homology-bockstein-pairing}
 \langle B\eta,s\rangle
       =\langle\eta,\omega\bmod n\rangle,
\end{equation}
where $B$ is the Bockstein for
$0\to\Z/n\xrightarrow{\times n}\Z/n^2\to\Z/n\to0$ and
$\omega\bmod n$ denotes its image in $H_2(J,\Z/n)$.
The two connecting maps use compatible sign conventions.
\end{lemma}

\begin{proof}
In the continuous bar complex with profinite integral coefficients,
lift a cycle representing $s$ to a chain $c$ and write $dc=nb$.
Then $b$ is a cycle representing $\omega$. Lift a cocycle representing
$\eta$ to a locally constant integral cochain $\widetilde\eta$.
The Bockstein is represented by $d\widetilde\eta/n$ modulo $n$, and
\[
 \left\langle\frac{d\widetilde\eta}{n},c\right\rangle
   =\left\langle\widetilde\eta,\frac{dc}{n}\right\rangle
   =\langle\widetilde\eta,b\rangle\pmod n.
\]
This is \eqref{eq:homology-bockstein-pairing}.
\end{proof}

Work on the common full-measure set $\Omega_{\mathrm{hom}}$ obtained
by intersecting $\Omega_+$ from \cref{lem:def-positive} with the
full-measure sets of \cref{prop:underlying-moment,lem:fab-primary-X}.
In particular, $X_\Gamma=1$ and every $(X^{\ab})_p$ is finite.
For $V$ inflated from $\Gamma$, smallness makes
$H^1(X\rt\Gamma,V)=\Hom_\Gamma(X,V)$ finite; \cref{prop:def-limit}
then gives finite-dimensional second cohomology.

\begin{proposition}
\label{prop:homology-cohomology-bounds}
On $\Omega_{\mathrm{hom}}$, for every $p\nmid|\Gamma|$ and every
irreducible $\F_p[\Gamma]$-module $V$, with trivial $X$-action,
\begin{equation}\label{eq:homology-annihilator-bound}
 \dim_{\F_p}H^2_{\omega_X}(X\rt\Gamma,V)
       -\dim_{\F_p}H^1(X\rt\Gamma,V)
 \le
       \sum_{E\in U}\dim_{\F_p}V^{\gamma_E}
                  -\dim_{\F_p}V^\Gamma.
\end{equation}
\end{proposition}

\begin{proof}
Fix $p$ and $V$; all dimensions in this proof are over $\F_p$.
Let $\mathcal T\subseteq H^2_{\omega_X}(X\rt\Gamma,V)$ be a
finite-dimensional subspace. As in the proof of \cref{prop:def-limit},
choose a finite level $\mathcal L$ containing the corresponding
extension groups and the split extension, and put
$\boldsymbol H=\boldsymbol X^{\mathcal L}$. There is a subspace
\[
 \mathcal T_H\subseteq H^2(H\rt\Gamma,V)^{\mathcal L}
\]
whose inflation identifies it with $\mathcal T$. The module $V$ is
eligible, and $C(\mathcal L,\boldsymbol H)$ has positive measure.
For the quotient map $\pi:X\onto H$, naturality gives
$\omega_H=\pi_*\omega_X$, so every class in $\mathcal T_H$
annihilates $\omega_H$. Set
\[
 H^2_{\omega_H}(H\rt\Gamma,V)^{\mathcal L}=H^2(H\rt\Gamma,V)^{\mathcal L} \cap H^2_{\omega_H}(H\rt\Gamma,V).
\]
Thus $\mathcal T_H\subseteq H^2_{\omega_H}(H\rt\Gamma,V)^{\mathcal L}$.

\medskip\noindent\textbf{If $p\nmid n$.}
Then $\omega_{H,p}=0$ and $ H^2_{\omega_H}(H\rt\Gamma,V)^{\mathcal L}=H^2(H\rt\Gamma,V)^{\mathcal L}$.
Taking logarithms in \eqref{eq:def-support-prime} gives
\begin{equation}\label{eq:homology-annihilator-finite-level}
 \dim H^2_{\omega_H}(H\rt\Gamma,V)^{\mathcal L}
 \le\dim H^1(H\rt\Gamma,V)
       +\sum_{E\in U}\dim V^{\gamma_E}-\dim V^\Gamma.
\end{equation}

\medskip\noindent\textbf{If $p\mid n$.}
For $\alpha\in H^2(H\rt\Gamma,V)^{\mathcal L}$ and
$\varphi\in\Hom_\Gamma(V,\Z/n)$, applying
\cref{lem:homology-bockstein-pairing} after restriction to $H$ gives
\begin{equation}\label{eq:homology-bockstein-annihilator}
 s_H\bigl(B(\varphi_*\alpha)\bigr)
   =\varphi\bigl(\langle\alpha|_H,\omega_H\bmod p\rangle\bigr).
\end{equation}
For irreducible $V$, the Bockstein conditions are thus equivalent to
$\alpha\in H^2_{\omega_H}(H\rt\Gamma,V)^{\mathcal L}$: if $V$ is nontrivial, both are
automatic, while if $V=\F_p$ is trivial, any nonzero $\varphi$
is injective. Hence $ H^2(H\rt\Gamma,V)^{\mathcal L,s_H}$ is the kernel of the cup-product map
\[
  H^2_{\omega_H}(H\rt\Gamma,V)^{\mathcal L}\longrightarrow H^1(H\rt\Gamma,V^\vee)^*,\quad
 \alpha\longmapsto\bigl(\beta\longmapsto s_H(\alpha\smile\beta)\bigr),
\]
where $(-)^*=\Hom_{\F_p}(-,\F_p)$ and the values are identified
with $\F_p$ as in \cref{subsec:positive-quotients}. Consequently,
\[
 \dim H^2_{\omega_H}(H\rt\Gamma,V)^{\mathcal L}
 \le\dim H^2(H\rt\Gamma,V)^{\mathcal L,s_H}+\dim H^1(H\rt\Gamma,V^\vee).
\]
Taking logarithms in \eqref{eq:def-support-dividing} bounds the
right side by
\[
 \dim H^1(H\rt\Gamma,V)
       +\sum_{E\in U}\dim V^{\gamma_E}-\dim V^\Gamma,
\]
proving \eqref{eq:homology-annihilator-finite-level} also in this case.

Finally, $\dim\mathcal T=\dim\mathcal T_H\le\dim H^2_{\omega_H}(H\rt\Gamma,V)^{\mathcal L}$,
and inflation injects $H^1(H\rt\Gamma,V)$ into
$H^1(X\rt\Gamma,V)$. Applying
\eqref{eq:homology-annihilator-finite-level} and varying $\mathcal T$
proves \eqref{eq:homology-annihilator-bound}.
\end{proof}

\subsection{Proof of the structural theorem}

\begin{proof}[Proof of \cref{thm:homology-structure}]
Fix a sample in $\Omega_{\mathrm{hom}}$ and a prime
$p\nmid|\Gamma|$, and put
\[
 \mathcal Q_{X,p}=H_2(X,\Z_p)/\langle\omega_{X,p}\rangle.
\]

\medskip\noindent\textbf{1. From cohomology to multiplicities.}
The coefficient sequence for multiplication by $p$ gives
\begin{equation}\label{eq:homology-uct}
 0\longrightarrow H_2(X,\Z_p)/pH_2(X,\Z_p)
 \longrightarrow H_2(X,\F_p)
 \longrightarrow (X^{\ab})_p[p]\longrightarrow0.
\end{equation}
For $V$ inflated from $\Gamma$, evaluation duality and coprime
averaging give
\[
 H^i(X\rt\Gamma,V)
       \simeq\Hom_\Gamma(H_i(X,\F_p),V)\qquad(i=1,2).
\]
These cohomology spaces are finite for every irreducible $V$.
Semisimplicity applied to the continuous dual of $H_2(X,\F_p)$
therefore shows that $H_2(X,\F_p)$ is finite.

Since $p\nmid|\Gamma|$, applying $\Hom_\Gamma(-,V)$ to
\eqref{eq:homology-uct} is exact. Restricting to classes annihilating
$\omega_X$ gives the short exact sequence
\begin{equation}\label{eq:homology-annihilator-exact}
 0\longrightarrow{}\Hom_\Gamma((X^{\ab})_p[p],V)
 \longrightarrow{}H^2_{\omega_X}(X\rt\Gamma,V)
 \longrightarrow{}\Hom_\Gamma(\mathcal Q_{X,p}/p\mathcal Q_{X,p},V)
       \longrightarrow0.
\end{equation}
Indeed, the last term consists precisely of the functionals on
$H_2(X,\Z_p)/pH_2(X,\Z_p)$ that vanish on the image of
$\omega_{X,p}$, since
\[
 \mathcal Q_{X,p}/p\mathcal Q_{X,p}
   =H_2(X,\Z_p)/(pH_2(X,\Z_p)+\langle\omega_{X,p}\rangle).
\]
The first term consists of classes vanishing on all of
$H_2(X,\Z_p)/pH_2(X,\Z_p)$, so these automatically annihilate
$\omega_X$.

The module $(X^{\ab})_p$ is finite by \cref{lem:fab-primary-X}.
Its multiplication-by-$p$ sequence shows that $(X^{\ab})_p[p]$ and
$(X^{\ab})_p/p(X^{\ab})_p$ have the same composition factors;
semisimplicity makes them isomorphic as $\F_p[\Gamma]$-modules.
Thus the first term of \eqref{eq:homology-annihilator-exact} has
dimension $\dim H^1(X\rt\Gamma,V)$. Taking dimensions gives
\begin{equation}\label{eq:homology-annihilator-difference}
 \dim_{\F_p}\Hom_\Gamma(\mathcal Q_{X,p}/p\mathcal Q_{X,p},V)
 =\dim_{\F_p}H^2_{\omega_X}(X\rt\Gamma,V)
              -\dim_{\F_p}H^1(X\rt\Gamma,V).
\end{equation}
On the other hand, the split augmentation sequence and Frobenius
reciprocity give
\begin{equation}\label{eq:homology-frobenius}
 \dim_{\F_p}\Hom_\Gamma(\overline{\mathcal K}_{U,p},V)
   =\sum_{E\in U}\dim_{\F_p}V^{\gamma_E}
            -\dim_{\F_p}V^\Gamma.
\end{equation}
By \cref{prop:homology-cohomology-bounds}, the Hom dimension in
\eqref{eq:homology-annihilator-difference} is at most that in
\eqref{eq:homology-frobenius}, for every irreducible $V$.
Dividing by $\dim_{\F_p}\End_\Gamma(V)$ compares their
multiplicities. Semisimplicity therefore gives an epimorphism
\[
 \overline\phi_{X,p}:\overline{\mathcal K}_{U,p}\onto_\Gamma
       \mathcal Q_{X,p}/p\mathcal Q_{X,p}.
\]

\medskip\noindent\textbf{2. Lifting and combining the primary components.}
Lift the images of a $\Z_p$-basis of $\mathcal K_{U,p}$ directly
through $H_2(X,\Z_p)\onto\mathcal Q_{X,p}/p\mathcal Q_{X,p}$,
and average the resulting continuous linear map over $\Gamma$.
This gives an equivariant map
\[
 \psi_{X,p}:\mathcal K_{U,p}\longrightarrow H_2(X,\Z_p).
\]
Its composite $\phi_{X,p}:\mathcal K_{U,p}\to\mathcal Q_{X,p}$
still reduces to $\overline\phi_{X,p}$ modulo $p$, since that map
was already equivariant. Its image is compact and hence closed;
its compact cokernel has zero reduction modulo $p$. Compact
Nakayama makes $\phi_{X,p}$ surjective. In particular,
\[
 \psi_{X,p}(\mathcal K_{U,p})+\langle\omega_{X,p}\rangle
       =H_2(X,\Z_p).
\]

The Chinese remainder theorem gives
\[
 \langle\omega_X\rangle=\prod_{p\mid n}\langle\omega_{X,p}\rangle,
 \qquad
 H_2(X,\widehat{\Z})/\langle\omega_X\rangle
       =\prod_{p\nmid|\Gamma|}\mathcal Q_{X,p}.
\]
Taking products of the $\phi_{X,p}$ proves
\eqref{eq:homology-structure-quotient}. Taking products of the
$\psi_{X,p}$ and adjoining the map $\Z/n\to H_2(X,\widehat{\Z})$
sending $1$ to $\omega_X$ proves \eqref{eq:homology-structure-plus}.
Its source is topologically finitely generated: $\mathcal K_U$ is
free of finite rank over $R_\Gamma$ and the additional
summand is cyclic. The same is therefore true of its quotient
$H_2(X,\widehat{\Z})$.
\end{proof}

\subsection{The vanishing consequence}

\begin{proof}[Proof of \cref{thm:homology-vanishing}]
Under the hypotheses of the theorem, the augmentation map in
\eqref{eq:homology-global-augmentation} is the identity of
$R_\Gamma$, so $\mathcal K_U=0$. Also $n=1$ implies $\omega_X=0$.
The epimorphism in \eqref{eq:homology-structure-quotient} then
forces $H_2(X,\widehat{\Z})=0$ almost surely.
\end{proof}

The corresponding question for arbitrary finite Galois extensions can
be stated as follows.

\begin{question}
\label{conj:arithmetic-virtual-pronilpotence-general}
Let $K/\mathbb Q$ be a finite Galois extension and put
$\Gamma=\operatorname{Gal}(K/\mathbb Q)$. Let $c\in\Gamma$ be
complex conjugation (so $c=1$ if $K$ is totally real). Assume that
\[
 \dim_{\mathbb C}W^c>0
 \quad\text{for every nontrivial irreducible complex $\Gamma$-module $W$.}
\]
Let $\pi$ be the set of primes not dividing $|\Gamma|$, and let $L/K$
be the maximal Galois pro-$\pi$ extension unramified at every finite
place and split at every real place. Does $\operatorname{Gal}(L/K)$ necessarily contain a pronilpotent open subgroup?
\end{question}

For general $\Gamma$ of odd order, we do not know whether the
analogous questions have affirmative answers even after replacing
``pronilpotent'' by ``prosolvable.''

\end{document}